\documentclass[11pt,oneside,reqno]{amsart}
\usepackage{epsfig,psfrag,subfigure}
\usepackage{amsmath,amssymb,amsfonts,latexsym, mathrsfs,euscript}
\usepackage{graphicx,graphics,cite}
\usepackage{wrapfig}
\usepackage[usenames]{color}
\usepackage{geometry}                
\usepackage{algorithm,algorithmicx,algpseudocode}

\usepackage{epstopdf}
\usepackage{subfig}
\usepackage{gensymb}

\DeclareGraphicsExtensions{.pdf,.png,.jpg}

\newtheorem{theorem}{\bf Theorem}

\newcommand{\beq}{\begin{equation}}
\newcommand{\eeq}{\end{equation}}
\newcommand{\beqa}{\begin{eqnarray}}

\newcommand{\eeqa}{\end{eqnarray}}

\newcounter{nfig}

\newcommand{\ignore}[1]{}

\newcommand\ds{\displaystyle}

\newcommand\bs{\boldsymbol}

\theoremstyle{remark}
\newtheorem{note}{Note}

\title[Geodesic landmark shooting algorithm]{A class of fast geodesic shooting algorithms for template matching and its applications via the $N$-particle system of the Euler-Poincar\'e equations}
\date{\today}

\author{Roberto Camassa}
\address{Department of Mathematics, University of North Carolina, Chapel Hill, 27599, USA}
\email{camassa@amath.unc.edu\\ TEL:919-962-8476}

\author{Dongyang Kuang}
\address{Department of Mathematics, University of Wyoming, Laramie, WY 82071-3036, USA}
\email{dkuang@uwyo.edu\\ TEL:307-766-4221}

\author{Long Lee}
\address{Department of Mathematics, University of Wyoming, Laramie, WY 82071-3036, USA}
\email{llee@uwyo.edu\\ TEL:307-766-4368}
\begin{document}
\maketitle

\begin{abstract}
The Euler-Poincar\'e (EP) equations describe the geodesic motion on the diffeomorphism group. For template matching (template deformation), the Euler-Lagrangian equation, arising from minimizing an energy function, falls into the Euler-Poincar\'e theory and can be recast into the EP equations. By casting the EP equations in the Lagrangian (or characteristics) form, we formulate the equations as a finite dimensional particle system. The evolution of this particle system describes the geodesic motion of landmark points on a Riemann manifold. In this paper we present a class of novel algorithms that take advantage of the structure of the particle system to achieve a fast matching process between the reference and the target templates. The strong suit of the proposed algorithms includes (1) the efficient feedback control iteration, which allows one to find the initial velocity field for driving the deformation from the reference template to the target one, (2) the use of the conical kernel in the particle system, which limits the interaction between particles and thus accelerates the convergence, and (3) the availability of the implementation of fast-multipole method for solving the particle system, which could reduce the computational cost from $O(N^2)$ to $O(N\log N)$, where $N$ is the number of particles. The convergence properties of the proposed algorithms are analyzed. Finally, we present several examples for both exact and inexact matchings, and numerically analyze the iterative process to illustrate the efficiency and the robustness of the proposed algorithms.

\end{abstract}

\begin{description}
\item[{\footnotesize\bf keywords}]
{\footnotesize Euler-Poincar\'e equations, geodesic motion, diffeomorphism group, template matching, template deformation, particle system, landmark points, feedback control iteration.}
\end{description}

\section{Introduction}

Research in template deformation has been prosperous in the last two decades. Starting from Grenander's deformable template models \cite{bib:Grenander93}, the study of Riemannian geometry of groups of diffeomorphisms and geodesics on manifolds have produced valuable theoretical results as well as actual algorithms in applications \cite{bib:Holm, bib:eqs, bib:bmty05, anatomy, bib:mty06, bib:dgm98,bib:mumford,bib:jm00}. Template matching, or template deformation,  is a common tool used in shape analysis. The many applications of shape analysis include image registration, pattern recognition, biomedical image analysis, morphometry, database retrieval, surveillance, biometrics, military target recognition and general computer vision \cite{image1, image2, bib:DM98, LR01,bib:sjml05}.

The Euler-Poincar\'e (EP) equations, also called the Euler equations for diffeomorphisms, are of general interest as evolution equations on Riemannian manifolds endowed with Sobolev metrics \cite{NP, bib:hs03}. Template matching can be formulated as finding the shortest or least expensive path of continuous deformation of one geometric object (reference template) into another one (target template). In this context, the time-dependent deformation process is the so-called geodesic evolution, and the derivation of the geodesic evolution equations falls into the Euler-Poincar\'e theory, which produces the EP equations \cite{anatomy, bib:hrty04,bib:Holm}. It is worth pointing out that despite the links between the EP equations and template deformation are well established \cite{anatomy, bib:bmty05, bib:hrty04,bib:my01,bib:jm00,bib:mumford}, the use of EP equations as numerical algorithms for template matching is only the beginning \cite{bib:MM2}, and has not been thoroughly investigated. 

Mathematically the EP equations describe geodesic motion on the diffeomorphism group, and are equivalent to the Euler-Lagrangian equations arising from minimizing an energy function defined for the deformation paths. A practical application for template matching is computational anatomy (CA) \cite{bib:HolmSingular1,anatomy}. For medical images in CA, instead of intensity on a mesh grid (pixel), geometry features of the medical images can be discretized into a set of the so-called {\it landmark points}. Therefore a template matching problem, in terms of landmark points, becomes a landmark-matching problem. That is, given two collections of points $X_1,\cdots, X_N$ and  $Y_1,\cdots, Y_N$, the matching problem is to find a time-dependent diffeomorphic path $\bs{\phi}(\cdot)$ that costs the minimum energy, such that $Y_k=\bs{\phi}(X_k)$, for $k=1,\cdots N$ \cite{bib:mumford, bib:hrty04}.  To link the landmark-matching process to the EP equations, we have shown that by casting the EP equations in the Lagrangian (or characteristics) form, i.e. following the characteristic variables, we can formulate EP equations as a finite-dimensional particle system of ordinary different equations (ODEs)  \cite{bib:ckl14}. This system of ODEs, referred to as the { $N$-particle finite-dimensional dynamical system}, or {\em $N$-particle system}, has two variables, the position variable, representing the positions of the particles, and the momentum variable, representing the momenta that drive the motion of the particles. Therefore for the full diffeomorphism group, a finite set of landmark points on the landmark space can be represented by the $N$-particles of the EP equations. The collections of the points $\{ X_1,\cdots, X_N \}$ and $\{ Y_1,\cdots, Y_N \}$ are the position variable of the $N$-particle system at two different times, respectively, and the landmark-matching problem is to find the {\em initial} momenta of the particles that drive the motion of the particles from one position $\{ X_1,\cdots, X_N \}$ to another $\{ Y_1,\cdots, Y_N \}$.  This idea of finding the proper initial momenta to match the evolution of landmark points of two different templates turns the usual initial value problem of the $N$-particle system into a conventional boundary value problem of template matching. The aim of this paper is to introduce a framework for designing a class of fast algorithms via the principle of feedback optimal control for template matching.

\section{Problem setting}

The formulation by Grenander et al. \cite{bib:Grenander93, bib:dgm98, bib:GM98} models the variations between two shapes or image templates by the action of Lie groups (diffeomorphisms) on manifolds. Template matching adopting this approach essentially considers an optimization problem with constraint:
\begin{equation}\label{eq:optimize}
\begin{split}
&\rho(I_0, I_1)^2 = \min\limits_{\bs{\phi}}\int_{0}^{1} E(t)\\
&\text{where}\, I_1={\bs{\phi}}({\bs u};I_0,1)
\end{split}
\end{equation}
In this setting, $I_0$ and $I_1$ are the point-distributions of the two shapes in a manifold $M$, where $I_0$ is called the reference template and $I_1$ is called the target template. $\int_{0}^{1} E(t)$ is a measure introduced in the tangent space $TM$, usually defined as the energy required for carrying $I_0$ to $I_1$ along a certain trajectory in $M$ \cite{bib:mumford}.  The geodesic map $\bs{\phi}({\bs u}; I_0, t)$ takes the initial value ${\bs u} \in T_{I_0}M$ as input carrying $I_0$ to another point in $M$ along the geodesic at time $t$ with respect to $E$. 
For exact image matching, the kinetic energy is defined by
\begin{equation}\label{eq:energy}
E(t)=\int_{X}E(x, t)d {\bs x}=\Vert {\bs u} \Vert^2_{\scriptscriptstyle\mathcal{L}}
\end{equation}
where ${\bs u}$ defines the velocity field for the geodesic flow ${\bs{\phi}}$, $\mathcal{L}$ is an inertia operator of the form of $\mathcal{L}=({\bs I}-\alpha^2 \Delta)^\nu$, where ${\bs I}$ is the identity matrix, $\Delta$ is the Laplacian, $\alpha$ is constant depending on the resolution level, and $\nu>0$ is a constant depending on the metric used, i.e. the order of Sobolev norm. The norm in tangent space is introduced as:
\beq\label{eq:u-square}
\Vert {\bs u} \Vert^2_{\scriptscriptstyle\mathcal{L}}=<{\bs u}, {\bs u}>_{\scriptscriptstyle{\mathcal{L}}}=<\mathcal{L}{\bs u}, {\bs u}>=\int_{D}\mathcal{L}{\bf u}\cdot{\bs u}\,\,d{\bs x},
\eeq
where $D$ is a compact subset, and 
\beq
{\bs u}(x,t) =\frac{\partial\bs{\phi}}{\partial t}\left(\bs{\phi}^{-1}(x, t), t\right).
\eeq

Traditionally, for images matching (intensity assigned to a mesh grid), or landmark matching (point distribution for the geometrical features of images), finding the minimum-energy deformation path is a boundary value problem, for which the image templates $I_0$ and $I_1$ are treated as two end points of a cylinder domain $M \times [0,1]$. Under this constraint, the optimization process updates the current curve (flow) with the two fixed end points $I_0$ and $I_1$, by using methods such as the steep (gradient) descent, to obtain the minimizer of the given energy function in Eq. (\ref{eq:optimize}),  The use of this diffeomorphisms (mappings) approach for analyzing images or shapes is mathematically sound and has been very successful \cite{bib:shooting12, bib:bmty05, anatomy, bib:jm00, bib:mty06}, but a major limitation here is the high computational cost \cite{bib:sjml05}.  

In this paper,  we introduce a class of geodesic shooting algorithms. The algorithms treat the template matching problem as an initial value problem with an unknown initial condition. Given a guess of the initial velocity (momenta), a sequence of approximate target templates, $I_1^{(k)}$, $k=0,1,\cdots$ are generated by solving the $N$-particle system of the EP equations. These approximate templates are used to obtain the corrections for the updates. Through this feedback-control loop, the initial velocities (momenta) that carry the landmark points of the reference template to that of the target template can be established.  

It is worth noting that the main differences between our algorithms and the traditional algorithms are: (1) philosophically, the previous methods search an approximate minimum-energy path to the exact target template, whereas in the proposed algorithms, the reference template follows an exact minimum-energy path to an approximate target template, driven by approximate initial momenta, (2) the proposed algorithms have the ability to predict the deformation beyond the target template, i.e. the algorithms could continue the deformation of the reference template beyond the target template by using the approximate initial momenta, and (3) since the kernels of the $N$-particle system have multipole expansions, the fast-multipole methods \cite{bib:hjz09} can be implemented for solving the particle system in the proposed algorithms.  

We also like to point out that an algorithm utilizing the $N$-particle system of the EP equations was introduced by McLachlan and Marsland for image registration \cite{bib:MM2}. For this algorithm, the initial conditions (momenta) that carry the landmarks from the reference template to target one are found by minimizing the Euclidean distance between the target template and the simulations, other than using the feedback-control iteration proposed in this paper. The minimization problem was solved by a sub-space trust region method based on the interior-reflective Newton method, provided by a MATLAB optimization package, and the metric of the $N$-particle system is the smooth Gaussian metric. Similar to the traditional minimization algorithms for template matching, this algorithm suffers from slow convergence and high computational cost, in particular for large numbers of landmarks \cite{bib:MM2}.

The major difference of the proposed algorithm from that of McLachlan and Marsland's is that our algorithm updates the initial guess of the velocity by using the vectors of the shooting error, until the error between the simulations and the target template is within a prescribed tolerance. This novel strategy reduces the computation cost by almost a half, compared with the approach of Newton's iteration, when a large number of landmarks are used, according to our numerical experiments.  Furthermore, a non-smooth kernel (the conical shape) that reduces the interaction between particles is used for the $N$-particle system to accelerate the convergence. Our results show that for the planar landmark matching problems, the proposed algorithms converge relatively fast without the implementation of a fast-multipole method for large numbers of landmarks, even for problems with sharp edges and problems with rotation and translation. If a fast-multipole algorithm is implemented, the computational cost could drop from $O(N^2)$ to at least $O(N\log N)$ in each time evolution \cite{bib:hjz09}, where $N$ is the number of landmark points (particles). 

\section{The EP equations and the Green's functions}

By using index notation with Einstein convention on sums over repeated indexes for 
the (column) vectors, $\bs{m}\equiv \left\{ m_{i} \right\}_{i=1}^n $ and $\bs{u}\equiv \left\{ u_{i} \right\}_{i=1}^n $, the EP equations can be written as
\begin{equation}\label{eq:EPDIFFi}
\partial_t m_{i} + u_{j}\partial_{j}m_{i}+ m_{j} \partial_{i} u_{j} + m_{i} \partial_{j} u_{j}=0 \, , 
\end{equation}
with $t\in\mathbb{R}^+$, $\bs{x}, \bs{u}$ and $\bs{m}\in\mathbb{R}^{n}$, and spatial partial derivatives are labeled by coordinate index. 
The velocity $\bs{u}$ and the momentum variable $\bs{m}$ are formally related by an elliptic  
 operator $\mathcal{L}$
\begin{equation}\label{eq:Yukawa}
\bs{m}=\mathcal{L}\bs{u} \, . 
\end{equation}
Let $\mathcal{L}$ be a self-adjoint operator $\mathcal{L}\equiv\mathcal{L}^{ \nu}$, where $\mathcal{L}^{\nu}$ is defined by
\begin{equation}\label{eq:elliptic}
\mathcal{L}^{\nu}=(\bs{I}-\alpha^2\nabla^2)^{\nu},
\end{equation}
parametrized by $\alpha^2$ and power $\nu >0$.  Here ${\bs I}$ is the identity matrix. For any $\nu>0$, including non-integer values, equation (\ref{eq:Yukawa}) is defined in the Fourier space 
\beq\label{eq:m-Fourier}
\begin{split}
& \bs{\hat{u}} = (\hat{\mathcal{L}^{\nu}})^{-1}\bs{\hat{m}},\quad\text{where}\,\,\,(\hat{\mathcal{L}^{\nu}})^{-1} = \frac{1}{(1+\alpha^2 |\bs{k}|^2)^{\nu}},\\ & |\bs{k}|=\sqrt{k_1^2+k_2^2\cdots+k_n^2},
\end{split}
\eeq
where $k_i$ is the  $i^{th}$ wavenumber. Since $\mathcal{L}^{ \nu}$ is rotationally invariant and diagonal, then $\bs{G}(\bs{x}) = G_{\nu-n/2}(|\bs{x}|)\bs{I}$ for a scalar function $G_{\nu-n/2}$, with $|\bs{x}|=\sqrt{x_1^2+x_2^2+\cdot+x_n^2}$, where $n$ is the dimension.  The scalar Green function $G_{\nu-n/2}$ is then written as
\beq\label{eq:Green-nD}
G_{\nu-n/2}(|\bs{x}|) = \frac{2^{n/2- \nu}}{(2\pi \alpha)^{n/2}\alpha^{\nu}\Gamma(\nu)}|\bs{x}|^{\nu -n/2}K_{\nu-n/2}\left(\frac{|\bs{x}|}{\alpha}\right),
\eeq
where $K_{\nu -n/2}$ is the modified Bessel function of the second kind of order $\nu -n/2$ and $\Gamma(\nu)$ is the usual notation for the Gamma function~\cite{bib:mumford}. A notable special parametric choice is the two-dimensional Green function's for $\alpha=1$ and $ \nu=3/2$, for which it takes the simple form 
\beq\label{eq:nu=3/2}
G_{1/2}(|\bs{x}|) = \frac{1}{2\pi}e^{-|\bs{x}|} \, . 
\eeq
We dub this Green's function as ``conon," due to the spatially conical shape of the function. In general, for $n=2$ (two-dimensional space), the regularity of the Green's function $2\pi G_{\nu-1}(r)$ is described as follows.

\begin{figure}[tbh]
\begin{center}
\includegraphics[width=4.0in]{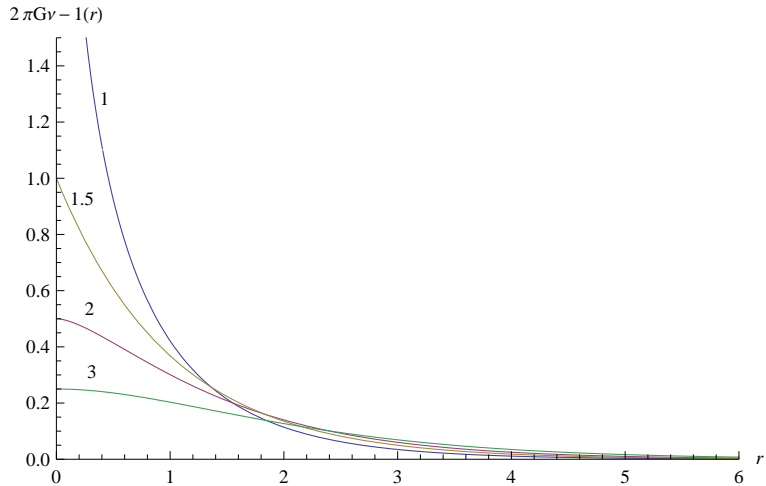}
\end{center}
\caption {Plots of $2\pi G_{\nu-1}(r)$ for $ b=1, 1.5, 2$ and $3$, where $G_{\nu-1}(r)$ is the two-dimensional Green's functions of the Yukawa operator $\mathcal{L}^{ \nu}$. $\alpha=1$ in the plots.}
\label{fig:Greens}
\end{figure}

\begin{enumerate} 
\item For the range $1/4< \nu\le 1$ the Green's function $G_{\nu-1}(|\bs{x}|)$ is unbounded. 

\item For the range $1< \nu <3/2$ the function is bounded but non-differentiable at the peak, with the radial derivative suffering an infinite jump there (cusp). 

\item At $\nu= 3/2$, the jump in radial derivative becomes finite. 

\item For the range $3/2 < \nu \le 2$ the derivative of the function is continuous, but with a infinite second derivative at the peak. 

\item Similar intervals can be defined for higher smoothness properties of the solution. In particular, for $2 < \nu<\infty$ the second derivative of the function is continuous. 
\end{enumerate}
The function $2\pi G_{\nu-1}(r)$ for the critical values $\nu=1, 1.5, 2,$ and 3, respectively is plotted in Fig. \ref{fig:Greens}.
Green's functions in the kernels of the $N$-particle system act like weight functions. From Figure \ref{fig:Greens}, we see that for a small $\nu$ value, the Green's function decays faster than that with a large $\nu$ value. This indicates that for a non-smooth Green's function, a particle has strong influence on other particles only when those particles are close enough to this particle. On the other hand, for a smooth Green's function, the influence of a particle to other particles spreads further away. Based on our numerical experiments, the choice of the metric $\mathcal{H}^{\nu}$ (corresponding to the Green's function  $G_{\nu-n/2}$) has an impact on the convergence rate of the geodesic shooting algorithms. In particular, we found that the metric $\mathcal{H}^{3/2}$ corresponding to the conical Green's function, $G_{1/2}$, provides the best convergence rate among the tested metrics. We remark that most of the previous work  in the literature, including \cite{bib:MM2}, used the smooth Gaussian kernel, $G_{\infty}$, corresponding to the metric $\mathcal{H}^{\infty}$ in their algorithms. We also note that metrics with $\nu<1.5$ are not suitable for our problems, since for those Green's functions, due to the low regularities, some kind of mollification will be needed for evaluating the kernels of the $N$-particle system.  We investigate the convergence property of the metrics $\mathcal{H}^{\nu}$ in Section \ref{sec:metric}.

\subsection{The EP equations and the $N$-particle system}

Previously we have shown that the EP equation (\ref{eq:EPDIFFi}) is the Eulerian counterpart of the following Lagrangian formulation \cite{bib:ckl14}:
\begin{equation}\label{eq:int-diff}
\begin{split}
{d \bs{q} \over d t} & = \int_{\mathbb{R}^{n}} 
G_{\nu-n/2}\big(|\bs{q}(\bs{\xi},t) -\bs{q}(\bs{\eta},t)|\big) 
\, \bs{p} (\bs{\eta},t) \, d V_\eta \, ,\\
{d \bs{p} \over d t} & = - \int_{\mathbb{R}^{n}} 
G'_{\nu-n/2}\big(|\bs{q}(\bs{\xi},t) -\bs{q}(\bs{\eta},t)|\big) 
{\bs{q}(\bs{\xi},t) -\bs{q}(\bs{\eta},t) \over |\bs{q}(\bs{\xi},t) -\bs{q}(\bs{\eta},t)|}
\bs{p}(\bs{\xi},t) \cdot \bs{p} (\bs{\eta},t) \, d V_\eta \,,
\end{split}
\end{equation}
where $\bs{q}$ is the position variable and $\bs{p}$ is the momentum variable. In this form, the equations of motion are a canonical Hamiltonian system with respect to variational derivatives $\delta/\delta\bs{q}$ and 
$\delta/\delta \bs{p}$ 
\begin{equation}
\dot{\bs{q}} (\bs{\xi},t) ={\delta H \over \delta \bs{p}}\, , 
\qquad \dot{\bs{p}} (\bs{\xi},t) =-{\delta H \over \delta \bs{q}}\,,
\label{eq:charsham}
\end{equation}
of the Hamiltonian functional 
\begin{equation}
H\equiv {1\over 2}\int_{\mathbb{R}^{n}} \int_{\mathbb{R}^{n}} 
G_{\nu-n/2}\big(|\bs{q}(\bs{\xi},t) -\bs{q}(\bs{\eta},t)|\big) \, \bs{p}(\bs{\xi},t) \cdot \bs{p}(\bs{\eta},t) \, 
dV_\xi \, dV_\eta  \, .
\label{eq:ham}
\end{equation}
Discretizing the above integral-differential equations and absorbing the grid sizes $dx$ and $dy$ into the momentum variable $p$, yields the $N$-particle system 
\beq\label{eq:N-particle}
\begin{split}
\ds\frac{d \bs{q}_i}{dt}& =\sum \limits_{j=1}^N G_{\nu-1}(|\bs{q}_i-\bs{q}_j|)\bs{p}_j,\\
\ds\frac{d\bs{p}_i}{dt}& = -\sum_{\substack{j=1\\ j\ne i}}^N(\bs{p}_i\cdot\bs{p}_j)G_{\nu-1}'(|\bs{q}_i-\bs{q}_j|)\frac{\bs{q}_i-\bs{q}_j}{|\bs{q}_i -\bs{q}_j|},
\end{split}
\eeq
and the discrete Hamiltonian
\beq\label{eq:H}
H = \frac{1}{2}\sum_{i=1}^{N}\sum_{j=1}^{N} (\bs{p}_i\cdot \bs{p}_j)G_{\nu-1}(\bs{q}_i-\bs{q}_j). 
\eeq
 Note that the velocity $\bs{u}(x, t)$ can be reconstructed at any time $t$ by 
\beq\label{eq:u}
\bs{u}(\bs{x}, t) = \sum_{j=1}^{N}G_{\nu-1}(|\bs{x} - \bs{q}_j(t)|)\bs{p}_j(t).
\eeq
From Eqs. (\ref{eq:energy}), (\ref{eq:u-square}), (\ref{eq:H}), and (\ref{eq:u}), we obtain
\beq\label{eq:E=2H}
 \begin{split}
E= &\int \mathcal{L}^{\nu} \sum\limits_{i=1}^N G_{\nu-1}(\bs{x}-\bs{q}_i)\bs{p}_i
       \cdot \sum\limits_{j=1}^N G_{\nu-1}(\bs{x}-\bs{q}_j)\bs{p}_j\ d{\bs x}\\
     = &  \sum_{i=1}^{N}\sum_{j=1}^{N} \int\mathcal{L}^{\nu}G_{\nu-1}(\bs{x}-\bs{q}_i)\bs{p}_i\cdot 
    G_{\nu-1}(\bs{x}-\bs{q}_j)\bs{p}_j\ d{\bs x}\\
    =& \sum_{i=1}^{N}\sum_{j=1}^{N} \int \bs{\delta}(\bs{x}-\bs{q}_i) {\bs{p}_i} \cdot G_{\nu-1}({\bs{x}-\bs{q}_j}){\bs{p}_j}d{\bs x}\\
    = & 2H,
 \end{split}
\eeq
where ${\bs{\delta}}$ is the Dirac delta function. It was shown that the function $\rho$ in Eq. (\ref{eq:optimize}) defines a metric \cite{anatomy}. 

Eq. (\ref{eq:E=2H}) implies that when $\bs{\phi}$ and $\bs{u}$ in Eq. (\ref{eq:optimize}) are determined, the Hamiltonian of the EP equations is $H=1/2 \rho^2$. In this paper, we use the Hamiltonian, corresponding to the velocity field that carries the flow from $I_0$ to $I_1$, as a measurement ($H$ is a semi-metric) that measures the difference between ${I}_0$ and  ${I}_1$ (or the logarithm of a conditional likelihood under the probability setting \cite{bib:dgm98}) in the tangent plane.

\section{Algorithm for Landmark Shooting}
We introduce a class of template matching algorithms in this section. We first consider the exact landmark matching problem 

\beq\label{eq:exact-matching}
I_1=\bs{\phi}({\bs u};I_0,\Delta t),
\eeq
where $I_0$ and $I_1$ are the shape or image templates represented by landmark points, and $\bs{\phi}$ is the geodesic flow determined by the corresponding Euler-Lagrangian equation.  After  the Euler-Lagrangian equation is recast in the Euler-Poincar\'e equation, the optimization problem becomes finding the initial condition for the initial value problem (\ref{eq:N-particle}). In Eq. (\ref{eq:exact-matching}), $\bs{u}$ is the initial vector field we intend to find, and $\Delta t$ is usually set to be 1 as a normalization. If the manifold $\bs{M}$ is compact, for given $I_1$, locally the existence and uniqueness of such $\bs{u}$ is guaranteed, and Eq. (\ref{eq:exact-matching}) is an one-to-one mapping \cite{bib:mumford}. 
The present algorithm is based on the idea of updating the velocity of the flow by using the difference between the landmarks of the target template and the solutions of the EP equations in each iteration. This process manifests the steep-decent method for the minimization problem. 

The iterative process is terminated at the $k^{th}$ iteration when $I_1^{(k)}$, calculated from the velocity ${\bs u}^{(k)}$,  satisfies $|I_1-I_{1}^{(k)}| < \epsilon$ for some prescribed threshold $\epsilon$. Here $I_1^{(k)}$ is a set, consist of the landmarks $\bs{q}_i^{(k)}$ that are solutions of the particle system of the EP equations. It is worth noting that since the new prediction is a linear combination of the previous one and a correction term from the observation quantity or pre-defined measurement, the iterative process of the algorithms is similar to the "state-space model" widely used in control theory. 
The geodesic landmark shooting algorithms are described in {\bf Algorithm 1}.

\begin{figure}[tbh]
\begin{center}
\includegraphics[width=4.0in]{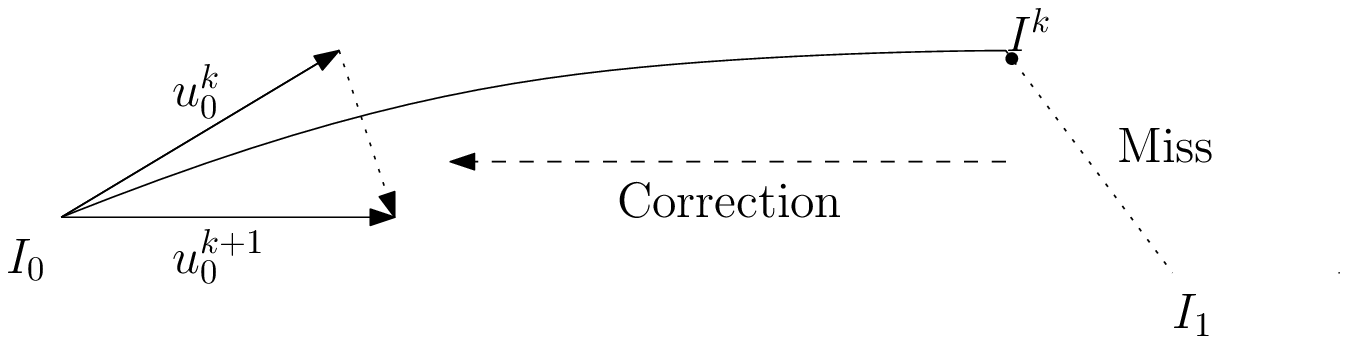}
\end{center}
\caption{Schematic iterpretation for {\bf Algorithm 1}.}
\end{figure}

\begin{algorithm}[t]
\begin{algorithmic}

\State Select corresponding landmarks (particles) between two shapes (images) $I_0$ and $I_1$;
\State Initialize the vector field $\bs{u}^{(0)}$ associated with the landmark points of $I_0$ with an initial guess of momenta $\bs{p}^{(0)}$;
\For{$k=0, 1, 2, \cdots$ until convergence}
\begin{enumerate}
\item[1.]  Evolve the $N$-particle system of EP equations (\ref{eq:N-particle}) to some fixed time ($t=1$), using the momenta associated with $\bs{u}^{(k)}$ to obtain the approximate landmarks $\bs{q}_{i}^{(k)}$ for  $I_1^{(k)}$, where $\bs{q}_i$ is the $i^{th}$ landmark, $i=1,\dots,N$.
\item[2.]  If the difference $||I_1-I_1^{(k)}|| < \epsilon$ for some norm, then stop.
\item[3.] The vector field is dynamically updated by $\bs{u}^{(k+1)} = \bs{u}^{(k)} + M^{(k)}\cdot(I_1-I_1^{(k)})$, where $M^{(k)}$ is a matrix.
\item[4.]  Update the momenta of the particle system $\bs{p}^{(k+1)}$ from  $\bs{u}^{(k+1)}$.
\end{enumerate}
\EndFor
\end{algorithmic}
\label{alg:proposed1}
\caption{Geodesic shooting algorithm for template matching}
\end{algorithm}


%
%
%
%
\begin{note}
\begin{enumerate}
\item There are many choices of placing the landmarks (selecting the position variable $\bs{q}$ of the particles in the particle system) onto the shape or image templates. For instance, we could place the particles on the entire grid points of a mesh, we could place the particles along the outline of a feature, we could position the particles only at the points-of-interest or significance in the image, such as edges and corners, or we could put the particles at places where the two images do not match.
\item The most ``efficient'' representation (i.e. a representation carries most of the information with the least number of landmarks.) varies, depending on different focus of the actual application. In this paper, for our numerical experiments we choose to place our landmarks along the contour around a given shape that mimics the outline of a substructure of an image. Potentially this approach could give us a clear picture about how biological structures deform.
  
\item The algorithm introduced here assumes that the chosen landmarks for the templates are appropriately labeled such that each ${\bs q}_{i} \in I_0$ will be carried to  ${\bs q}_{i} \in I_1$ by the flow when the matching process is complete. In practice, the biological structure will aid the labeling, but this is out of the scope of the current paper, and we will not pursue further here. It is worth noting that there are algorithms developed for unlabeled landmark matching. We refer readers to  the literature \cite{bib:unlabel1, bib:unlabel2} for more information.
  
\item Our initial guess for the velocity field is normally $\bs{u}^{(0)}=0$.  This is a reasonable initialization to test the stability of our algorithm. That is, if the problem is viewed under a stochastic setting: $I_0$ and $I_1$ is two realizations of certain process, their difference will be a zero mean process, $E(I_0-I_1)=0$. Thus, zero is an unbiased estimation to begin with.

\item Another way to find the velocity vector field is to ``directly'' solve the non-linear equation $I_1-\bs{\phi}({\bs u};I_0,\Delta t)=0$ for $\bs{u}$, by some iterative methods, such as the Newton iteration. Comparison in our numerical experiments shows that our choice of the optimizer is advantageous over the ``direct'' approach.

\item The present algorithm involves the evaluation of a double summations in the $N$-particle system. Suppose that the system (\ref{eq:N-particle}) is solved by an explicit time integrator, such as the classical Runge-Kutta algorithm, in which the summations are evaluated directly, the number of operations is $O(N^2)$, where $N$ is the number of  landmark points. This makes the shooting algorithms undesirable for dealing with a large database of shapes or images .The fast-multipole-method (FMM) developed by Huang et al. \cite{bib:hjz09} for the screened Coulomb interactions of $N$ particles in principle can be modified to solve the $N$-particle system in {\bf Algorithm 1}. The FMM in principle will reduce the computational cost to at least $O(N\log N)$. 
\end{enumerate}

\end{note} 
  
%
%
%
%
%

\section{Local convergence theorem}\label{sec:convergence}
%
%
%
%

In this section we present a local convergence theorem. Based on this theorem, we select the matrix $M^{(k)}$ in {\bf Algorithm 1} for our numerical experiments.
\begin{theorem}\label{them:1}
Let $\bs{\phi}:\mathbb{R}^{d}\mapsto\mathbb{R}^{d}$ be a diffeomorphism. Suppose that the first-order Taylor series expansion of $\bs{\phi}$ is uniformly bounded in the sense that 
\beq\label{eq:cond}
\begin{split}
{\bs{\phi}}(v)-{\bs{\phi}}(V) = & D\bs{\phi}(v)\cdot (v-V) + R(v-V),\\
\text{where}\,\,&||R(v-V))||\leq B||v-V||,
\end{split}
\eeq
for some constant $K\in\mathbb{R}$ in a region $\Omega := \{v, V\in \mathbb{R}^{d}:||v-V||<r\}$ for some constant $r\in\mathbb{R}$. Here $R(\cdot)$ is the remainder of the Taylor series and $D\bs{\phi}$ is the Jacobian matrix of the flow ${\bs{\phi}}$. Suppose that $D\bs{\phi}$ is non-singular. Then for the velocity update $\bs{u}^{(k+1)}=\bs{u}^{(k)}+M^{(k)}\cdot(I_1-I_1^{(k)})$, there exist a sufficiently large number  $N$, for which the error $||e_k||=||\bs{u}-\bs{u}^{(k)}|| \leq \epsilon$ for $k> N$ and $\epsilon<\!\!< r$. Here $\bs{u}$ is the true initial velocity for the reference template, and $\bs{u}^{(k)}$ is the $k^{th}$ iteration of the velocity vector found by  \bf{Algorithm 1}.
\end{theorem}

\noindent
{\em  Proof}: Let $\bs{I}$ be the identity matrix. Consider the $(k+1)^{th}$ update for the velocity field: $\bs{u}^{(k+1)}=\bs{u}^{(k)}+M^{(k)}\cdot(I_1-I_1^{(k)})$ in {\bf Algorithm 1}, where ${\bs{\phi}}(\bs{u};I_0,1)=I_1$, and $I_1^{(k)}$ is the solution of the EP equations in the $k^{th}$ interation.  Here $I_0$, $I_1$, and  $I_1^{(i)}$ are vectors of landmark points.  Suppose that condition (\ref{eq:cond}) is satisfied, the bounded first-order Taylor  expansion for ${\bs{\phi}}$ on the convex open domain $\{v, V\in \mathbb{R}^{d}:||v-V||<r\}$ for some constant $r\in\mathbb{R}$ at the $k^{th}$ iteration yields the following calculations:
\begin{eqnarray}
&\bs{u}^{(k+1)}-\bs{u}^{(k)}=M^{(k)}\cdot\left(I_1-I_1^{(k)}\right)\\
\Rightarrow  &\left(\bs{u}^{(k+1)}-\bs{u}\right)-\left(\bs{u}^{(k)}-\bs{u}\right)=M^{(k)}\cdot\left(\bs{\phi}(\bs{u})-{\bs{\phi}}(\bs{u}^{(k)})\right)\\
\Rightarrow  & e_{k+1}-e_k =-M^{(k)}\cdot D\bs{\phi}(\bs{u}^{(k)})\cdot e_k + M^{(k)}\cdot R(e_k)\\
\Rightarrow  & e_{k+1}=\left(\bs{I}-M^{(k)}\cdot D\bs{\phi}(\bs{u}^{(k)})\right)\cdot e_k + M^{(k)}\cdot R(e_k)\\ \label{c_rate}
\Rightarrow  & ||e_{k+1}||\leq||\bs{I}-M^{(k)}\cdot D\bs{\phi}(\bs{u}^{(k)})|| \,||e_k||+ B ||M^{(k)}||\,||e_k||\\
\Rightarrow  & \displaystyle\frac{||e_{k+1}||}{||e_k||}\leq||\bs{I}-M^{(k)}\cdot D\bs{\phi}(\bs{u}^{(k)}) || +B ||M^{(k)}||\label{eq:ratio}
\end{eqnarray}
Since $D\bs{\phi}$ is non-singular, $(D\bs{\phi})^{-1}$ exists. Choose $M^{(k)}=h(D\bs{\phi})^{-1}(\bs{u}^{(k)})$, where $0<h<1$ and assume that $||(D\bs{\phi})^{-1}(\bs{u}^{(k)})||\le c_k$, for some  $c_k$. Let  $C=\underset{k}{\max}\,\,c_k$, then  Eq. (\ref{eq:ratio}) becomes
\begin{equation}
\frac{||e_{k+1}||}{||e_k||}\le||\bs{I}-h\bs{I}|| +hCB \le |1-h| +hCB.
\end{equation} 
Thus there exists some $h$ so that $||e_{k+1}|| < \beta ||e_{k}||$ for $\beta= |1-h| +hCB<1$, provided $CB<1$. Therefore,  there exists a sufficiently large number $N$ such that $||e_k|| \le \epsilon$ for $k>N$, where $\epsilon <\!\!< r$.
\begin{note}
\begin{enumerate}
\item Theorem \ref{them:1} suggests that if the iteration $\bs{u}^{(k)}$ is close to the true velocity $\bs{u}$ in some norm at all times, then $\bs{u}^{(k)}$ converges to $\bs{u}$ for $k\geq N$ with $N$ sufficiently large, provided $M^{(k)}= h(D\bs{\phi})^{-1}(\bs{u}^{(k)})$ is bounded at all times.
\item In practice, it is too cumbersome and time-consuming to numerically compute $(D\bs{\phi})^{-1}(\bs{u}^{(k)})$. Instead, we let $M^{(k)}$ be a constant matrix in the form of $M^{(k)}=M=h{\bs I}$ with $0<h<1$ for all $k$ (analogous to assuming that $D\bs{\phi}(\bs{u}^{(k)}) = (D\bs{\phi})^{-1}(\bs{u}^{(k)}) =\bs{I}$, where $\bs{I}$ is the identity matrix), then it implies $C=1$, and therefore $CB<1$ is possible for a sufficiently small $B$. Indeed, our numerical convergence study shows that it suffices to use the constant matrix $M$ for all our numerical experiments. 
\item See Section \ref{sec:app} for the derivation of an theoretic optimal $M$ from the view point of minimizing variance estimation.
\end{enumerate}
\end{note}

\section{The inexact landmark matching}
Deviated from the exact template matching, the so-called inexact template matching \cite{anatomy, bib:mty06, bib:hrty04,bib:Holm} adds a running norm squared term to the energy functional (\ref{eq:energy}). The resulting functional, under the setting of landmark evolution, is written in the form of
\begin{equation}\label{eq:inexact-E}
\tilde{E} = \int_{0}^{1} \Vert \bs{u}^2 \Vert_{\mathcal{L}} dt + 
\frac{1}{\sigma^2}\sum \limits_{j=1}^N \Vert I_{1,j}-\bs{\phi}(\bs{u};I_{0,j},1)\Vert^2_{\mathbb{R}^d},
\end{equation}
where $I_{k,j}$ represents the $j^{th}$ landmark in template $k$. The minimizer of Eq. (\ref{eq:inexact-E}) provides the velocity to evolve the reference template $I_0$ to the target template $I_1$. This process of inexact matching is referred to as a metamorphosis, which provides a mechanism that allows the evolution to deviate from pure geodesic deformations \cite{bib:ty05, bib:Holm, bib:eqs}. The finite-dimensional particle system corresponding to the Euler-Lagrangian equations, arising from the variational problem of Eq. (\ref{eq:inexact-E}), can be written as follows \cite{bib:Holm, bib:eqs}:
\beq\label{eq:N-particle-inexact}
\begin{split}
\ds\frac{d \bs{q}_i}{dt}& =\sum \limits_{j=1}^N G_{\nu-1}(|\bs{q}_i-\bs{q}_j|)\bs{p}_j +\sigma^2 {\bs{p}_i},\\
\ds\frac{d\bs{p}_i}{dt}& = -\sum_{\substack{j=1\\ j\ne i}}^N(\bs{p}_i\cdot\bs{p}_j)G_{\nu-1}'(|\bs{q}_i-\bs{q}_j|)\frac{\bs{q}_i-\bs{q}_j}{|\bs{q}_i -\bs{q}_j|}.
\end{split}
\eeq
The parameter $\sigma^2$ provides a slightly inexact advection and the magnitude of  $\sigma^2$ is a weight defining how close to a pure geodesic deformation a metamorphosis will be \cite{bib:eqs}. The algorithm for inexact matching is similar to {\bf Algorithm 1}, except the stopping criterion used for inexact matching is the difference of approximate initial momenta computed between the $k^{th}$ and $(k-1)^{th}$ iterations. If the difference is smaller than the prescribed tolerance in some norm, then we consider the process is numerically convergent and stop the iteration. 


\section{Numerical Experiments - Exact Matching}\label{sec:exact}

In this section, we present several numerical experiments and study their numerical convergence properties to illustrate the advantages of {\bf Algorithm 1}. In all the experiments, unless specified otherwise, 64 landmarks ($N=64$) are used to represent the outline of a given shape.  We use the tolerance $\epsilon=10^{-3}$ as our stopping criterion. The metric of the self-adjoint operator $\mathcal{L}$ is $\mathcal{H}^{3/2}$ with $\alpha =1$ (i.e. $\mathcal{L}^{3/2} = (\bs{I} - \nabla^2)^{3/2}$), for which the Green's function is normalized and is written as  $G_{3/2}(|\bs{x}|)=e^{-|\bs{x}|}$. We use the Hamiltonian (\ref{eq:H}) to measure the potentials between the reference template and the target one. This semi-metric is invariant under a rotation, a translation, but not a scaling, according to our experiments (e.g. see Table \ref{tab:isoH}). We refer readers to the literature \cite{anatomy} for other invariant metrics and their properties. Finally, the numerical examples are computed by using Matlab(2011a) on Windows 7 Home Premium (SP1) platform with Intel\textsuperscript{\textregistered} i5-2430 processor @ 2.40GHz.

%
 
\subsection{Example 1} 
We construct the diffeomorphism between a circle and a rotated-and-shifted ellipse by using the proposed {\bf Algorithm 1} in this example. The reference template is a circle centered at origin with radius 2, and the target template is an ellipse with $a=4, b=1$ rotated by $45^{\degree}$ and shifted in positive x direction by 1 unit. This example illustrates the flexibility of {\bf Algorithm 1} for handling deformations with translations and rotations. The reference template $I_0$ and the target template $I_1$ are described as follows.
\beq\label{eq:I0} 
I_0: x^2+y^2=2^2,
\eeq
\beq\label{eq:ellips}
I_1:
\begin{cases}
x =& 4\cos(\frac{\pi}{4})\cos(\theta)+\sin(\frac{\pi}{4})\sin(\theta)+1, \\
y =& \cos(\frac{\pi}{4})\sin(\theta)-4\sin(\frac{\pi}{4}).
\end{cases}
\eeq
\begin{figure}[h]
\subfigure[t=0, $H=0$]{\includegraphics[width=2.in]{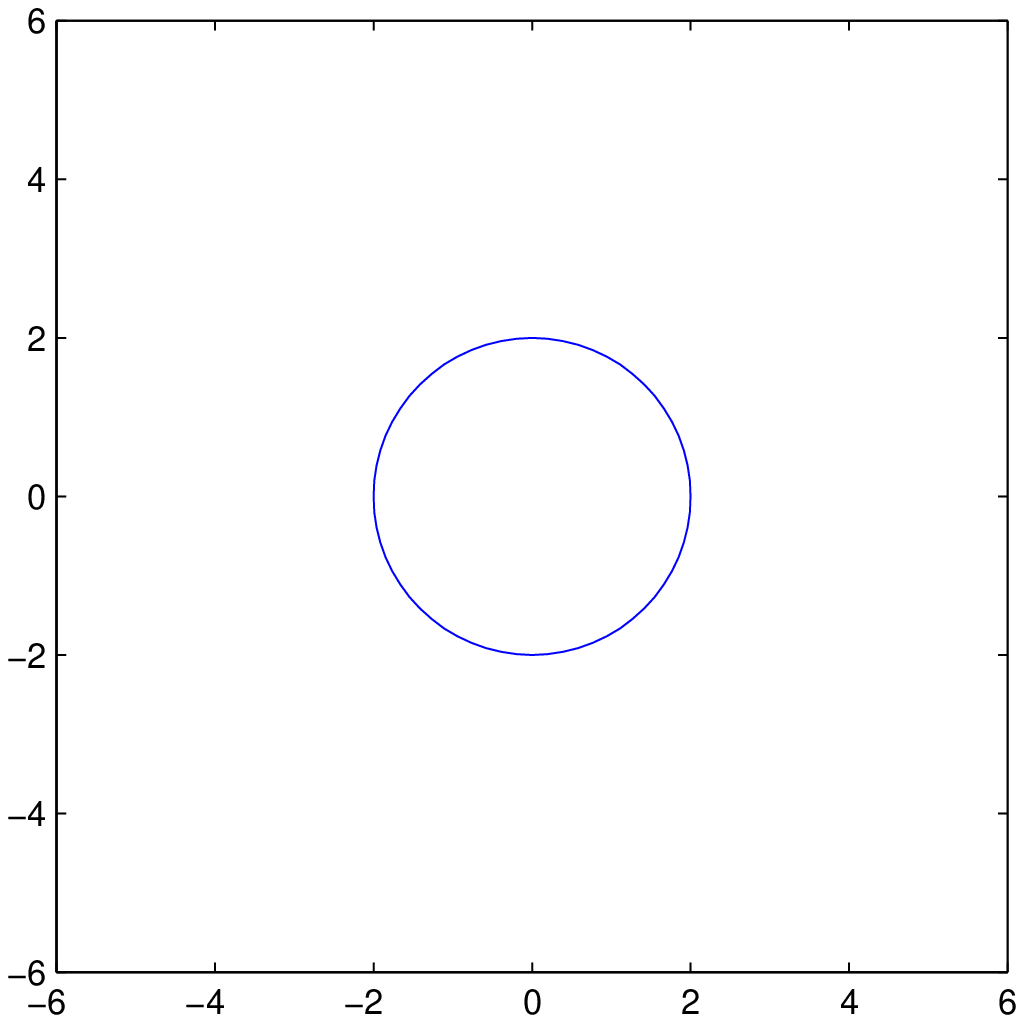}}
\subfigure[t=0.2, $H=1.8595$]{\includegraphics[width=2.in]{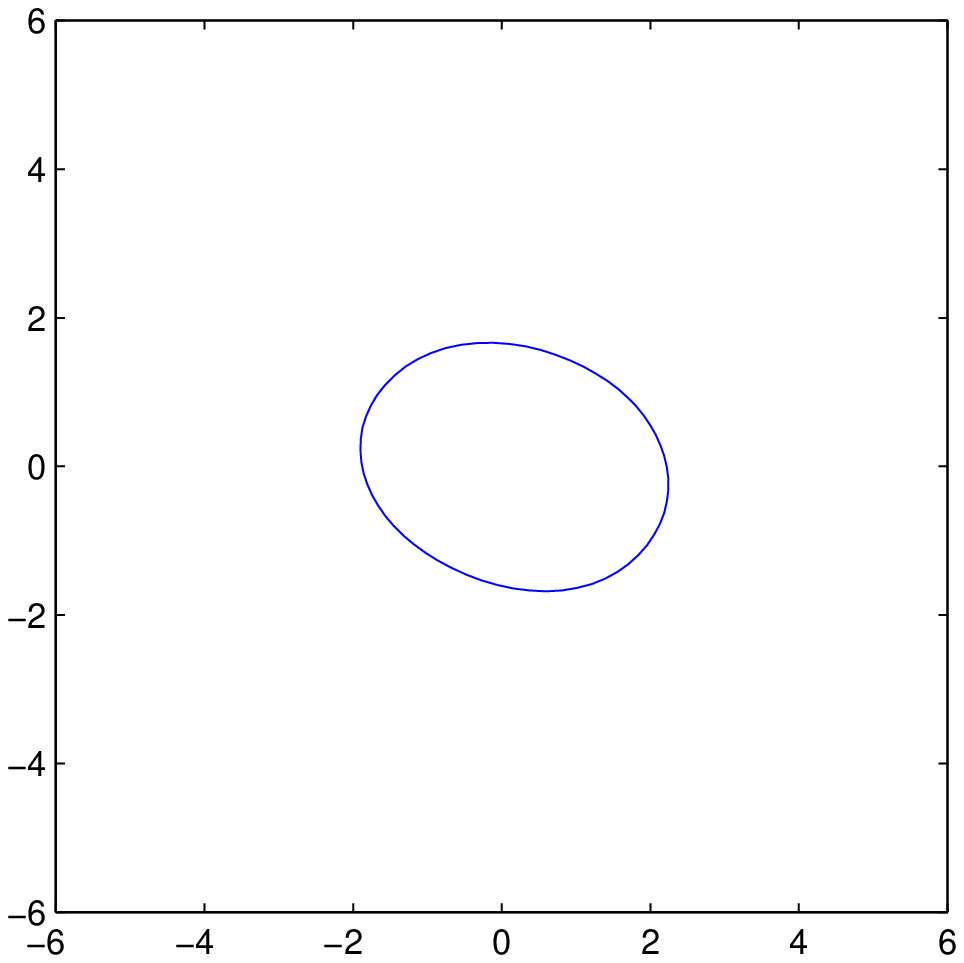}}
\subfigure[t=0.4, $H=7.4417$]{\includegraphics[width=2.in]{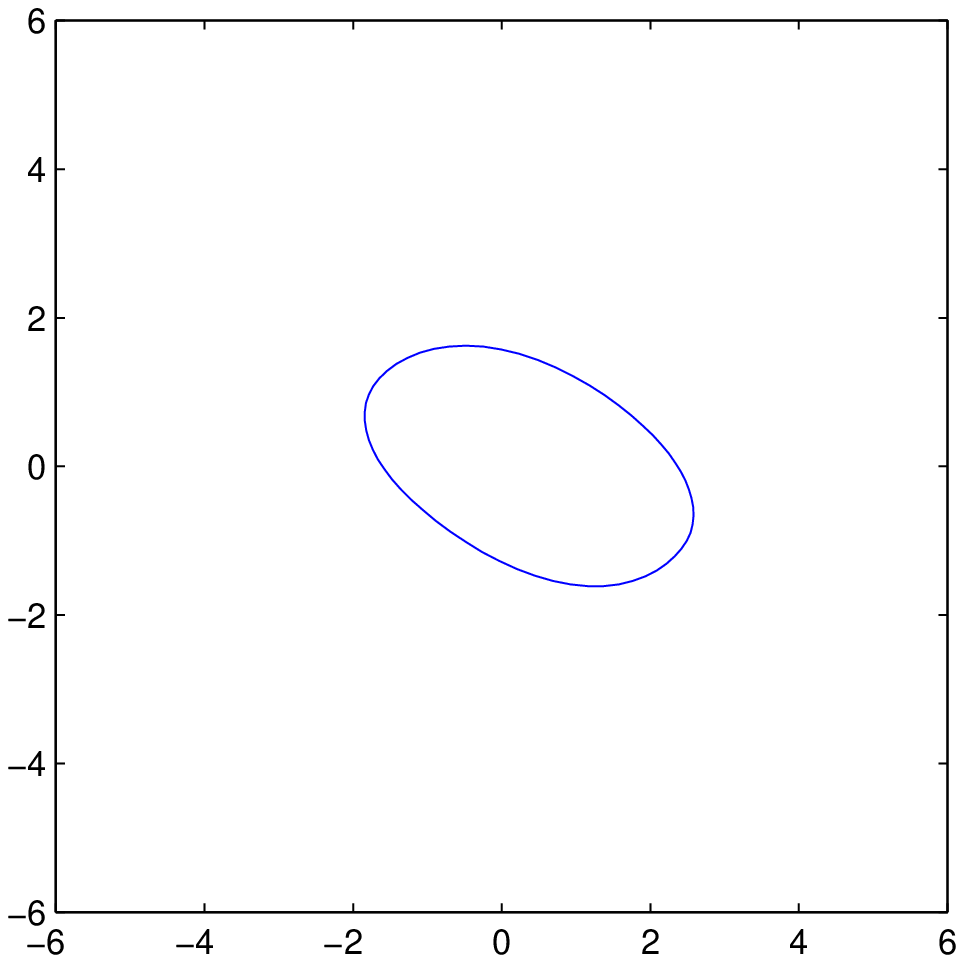}}\\
\subfigure[t=0.6, $H=16.7461$]{\includegraphics[width=2.in]{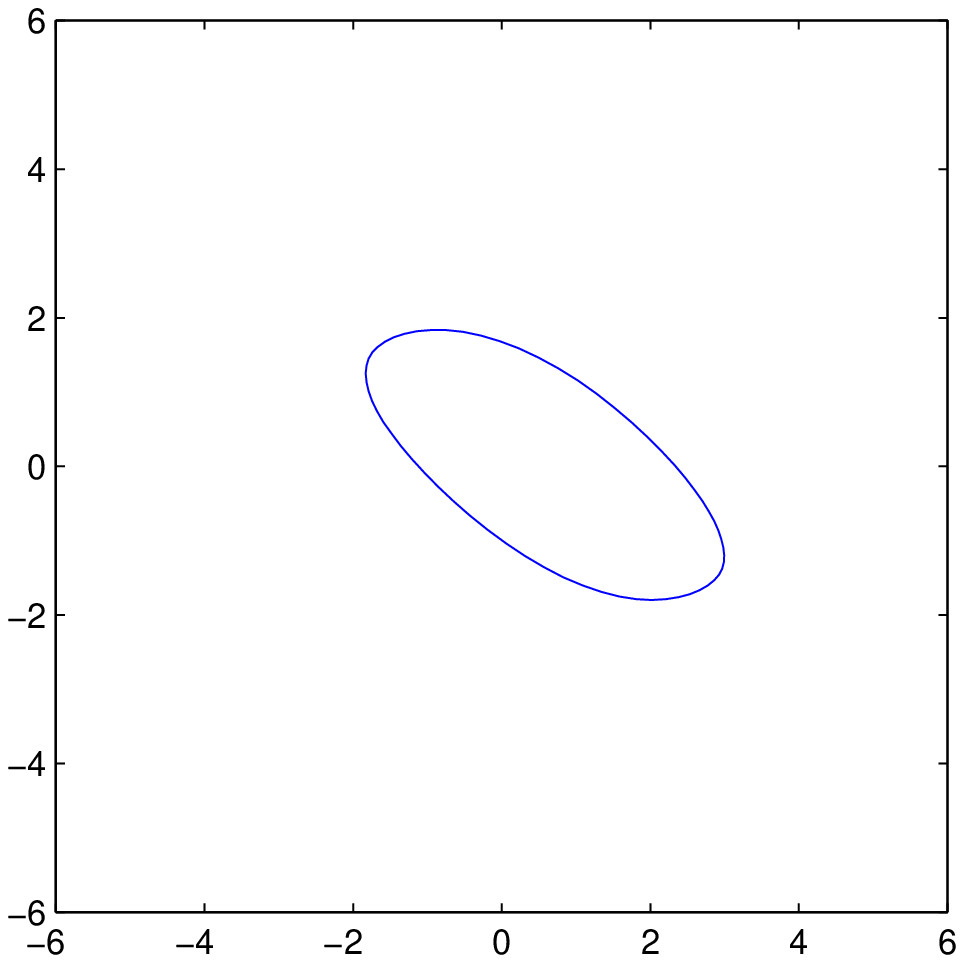}}
\subfigure[t=0.8, $H=29.7594$]{\includegraphics[width=2.in]{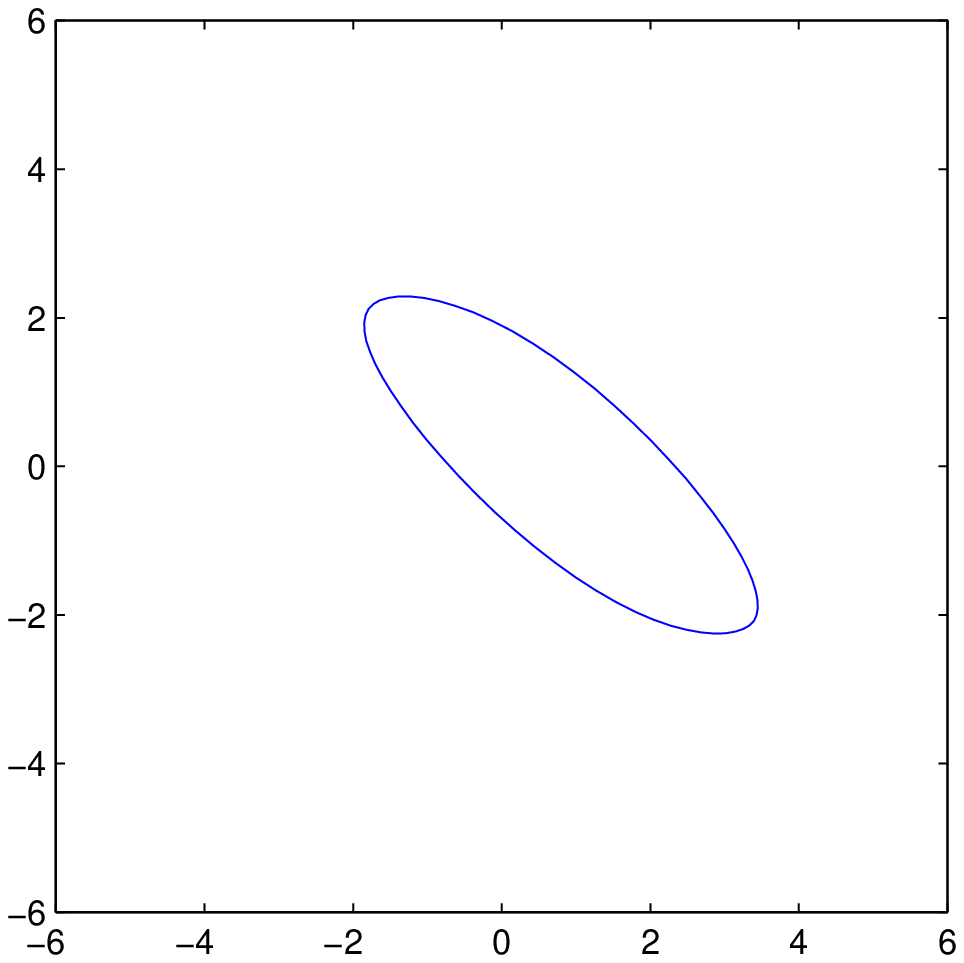}}
\subfigure[t=1, $H=46.5022$]{\includegraphics[width=2.in]{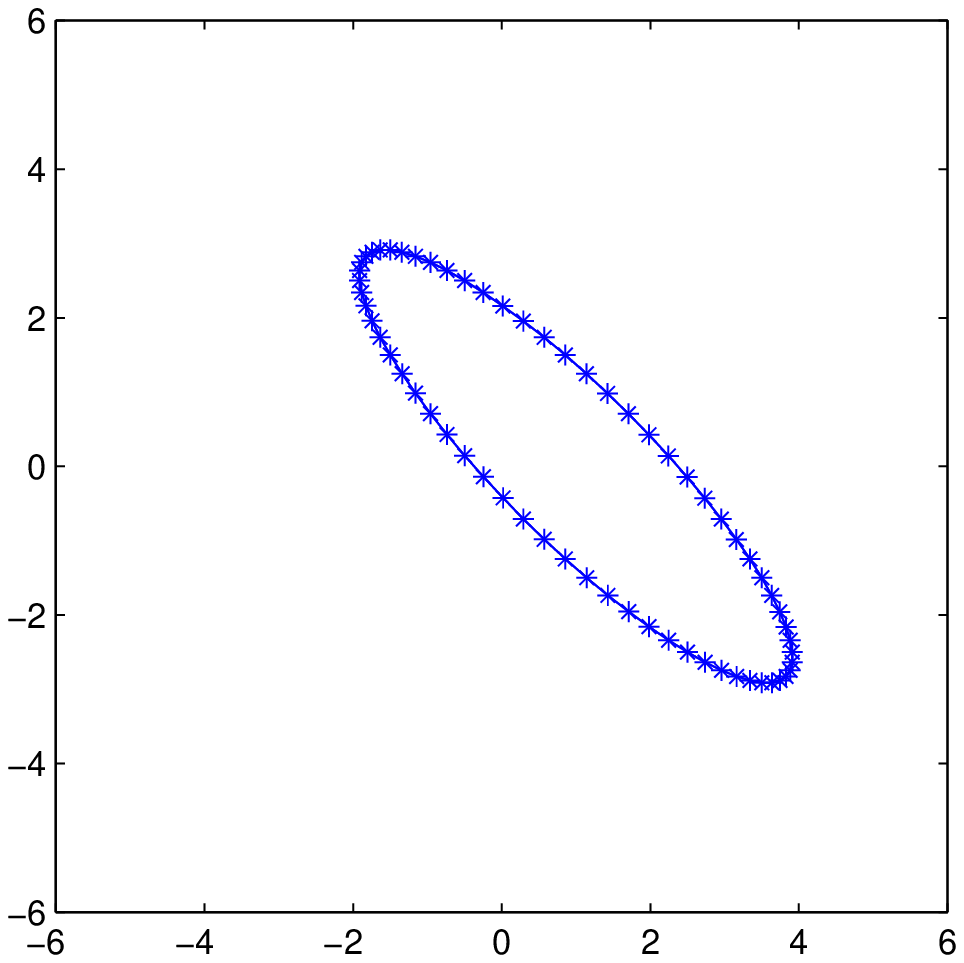}}
\caption{Diffeomorphism from a circle to a shifted and rotated ellipse. (a) is the reference template. (b)-(e) are the interpolative shapes along the path of minimum energy found by {\bf Algorithm 1}. The dotted line in (f) is the landmark distribution of the target template, whereas the solid line is the approximate target template.}
\label{ellipse}
\end{figure}

The diffeomorphism from $I_0$ to $I_1$ goes through rotation, dilation and translation. Figure (a) is the reference template. Figure \ref{ellipse} (b)-(e) show the interpolative stages during the evolution. The dotted line in Figure (f) is the landmark distribution of the target template, whereas the solid line is the approximate target template. $M=h\bs{I}$, where $h=0.3$ is used. It takes 197 iterations to reach the tolerance. 


To demonstrate the efficiency of {\bf Algorithm 1}, Table \ref{tab:comp1} compares {\bf Algorithm 1} with the approach of solving the non-linear equation $I_1-\bs{\phi}({\bs u};I_0,\Delta t)=0$ for $\bs{u}$ directly by the MATLAB intrinsic function {\bf fsolve}, a trust-region Newton algorithm. Table \ref{tab:comp1}(a) is the errors and elapsed CPU times when the stop criterion is achieved for various $N$, the number of landmarks, by using {\bf fsolve}, while Table \ref{tab:comp1}(b) shows the results by  {\bf Algorithm 1}. From the comparison, we see that when $N$ is small, the two methods are comparable, but when $N$ is large, {\bf Algorithm 1} is advantageous. The default stopping criterion of {\bf fsolve} is used for Table \ref{tab:comp1}(a). The stopping criterion for Table \ref{tab:comp1}(b) is the corresponding $\ell_2$ error found in Table \ref{tab:comp1}(a) for each $N$. The table also shows that the computational cost for {\bf Algorithm 1} is roughly  $O(N^2)$.


\begin{table}[h]
\caption{Comparison between (a) MATLAB {\bf fsolve}  and (b) {\bf Algorithms 1} with the search length $h=0.3$ for the update matrix $M=h\bs{I}$.}
\label{tab:comp1}
\begin{minipage}{0.45\textwidth}
\centering
(a)
\begin{tabular}{|l||c|c|}
\hline
\hline
$N$ & error(s) & Total time(s)\\
\hline
30 &  1.2006e-8 & 11.4370\\ 
\hline
40 &  5.1522e-10   & 26.2334\\
\hline
50 &  2.0985e-7 & 41.3252\\
\hline
60 & 1.1505e-7 & 77.0952\\
\hline

\hline

\end{tabular}
\end{minipage}
\begin{minipage}{0.45\textwidth}
\centering
(b)
\begin{tabular}{|l||c|c|}
\hline
\hline
$N$ & error(s) & Total time(s)\\
\hline
30 &  9.9032e-9 & 12.9654\\ 
\hline
40 &  4.5708e-10   & 28.5737\\
\hline
50 &  1.9258e-7 & 29.9143\\
\hline
60 & 8.8837e-8 & 43.8872\\
\hline

\hline

\end{tabular}
\end{minipage}
\end{table}

\subsection{Example 2} 
In the second example, we compute the diffeomorphism between a circle and a heart curve and illustrate some numerical properties of {\bf Algorithm 1}.  Similar to {\bf Example 1}, the reference template is a circle defined by
\beq\label{eq:I0-repeat} 
I_0: x^2+y^2=2^2,
\eeq
whereas the deformed target template is a $4^{th}$-order heart curve defined as
\beq \label{eq:heart}
 I_1: 
 \begin{cases}
 x= &\frac{1}{5}(13\cos(\theta)-5\cos(2\theta)-2\cos(3\theta)-\cos(4\theta))\\
 y= &\frac{1}{5}(16\sin(\theta)^3).
 \end{cases}
\eeq
The diffeomorphism involves the evolution of a convex region to a non-convex one, and the formation of a singularity from a piece of smooth curve. This example demonstrates the ability of the proposed algorithm to construct interpolative shapes between smooth (reference template) and sharp-edge (target template) planar curves.
\begin{figure}[tbh]
\subfigure[$t=0$, $H=0$]{\includegraphics[width = 2in]{fig1.eps}} 
\subfigure[$t=0.2$, $H=0.7404$]{\includegraphics[width = 2in]{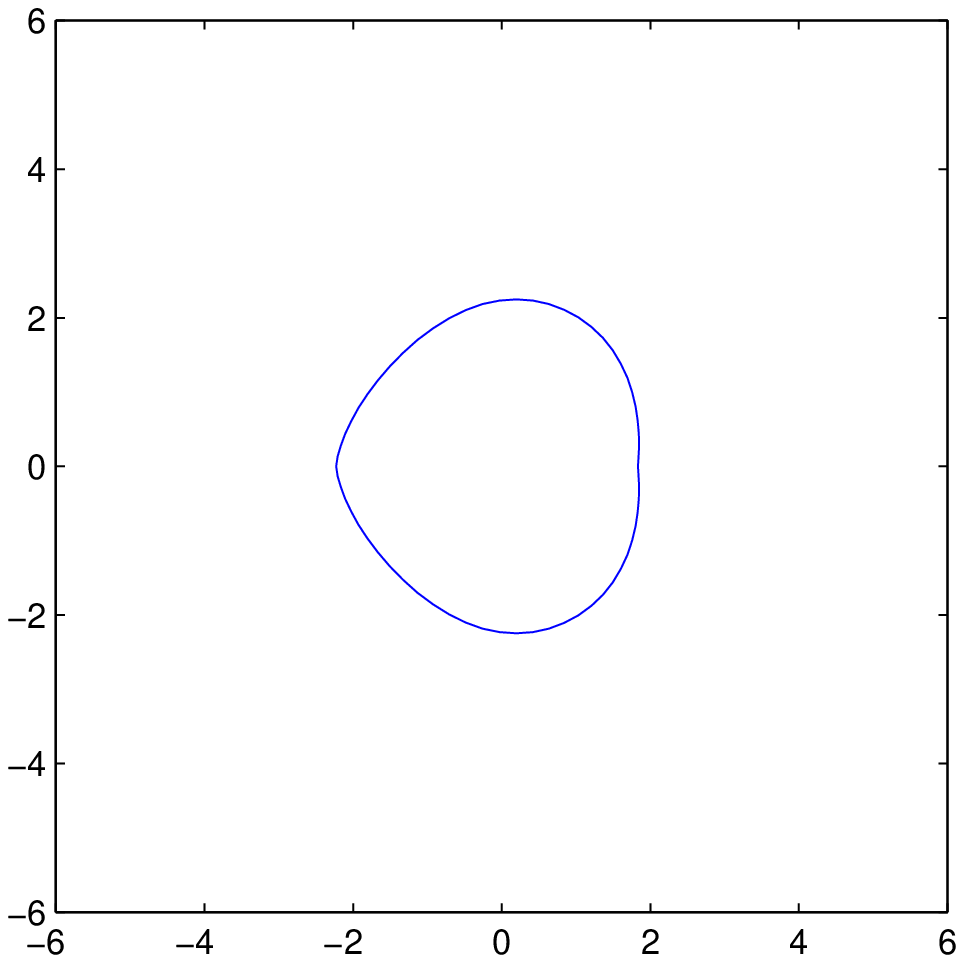}}
\subfigure[$t=0.4$, $H=2.9693$]{\includegraphics[width = 2in]{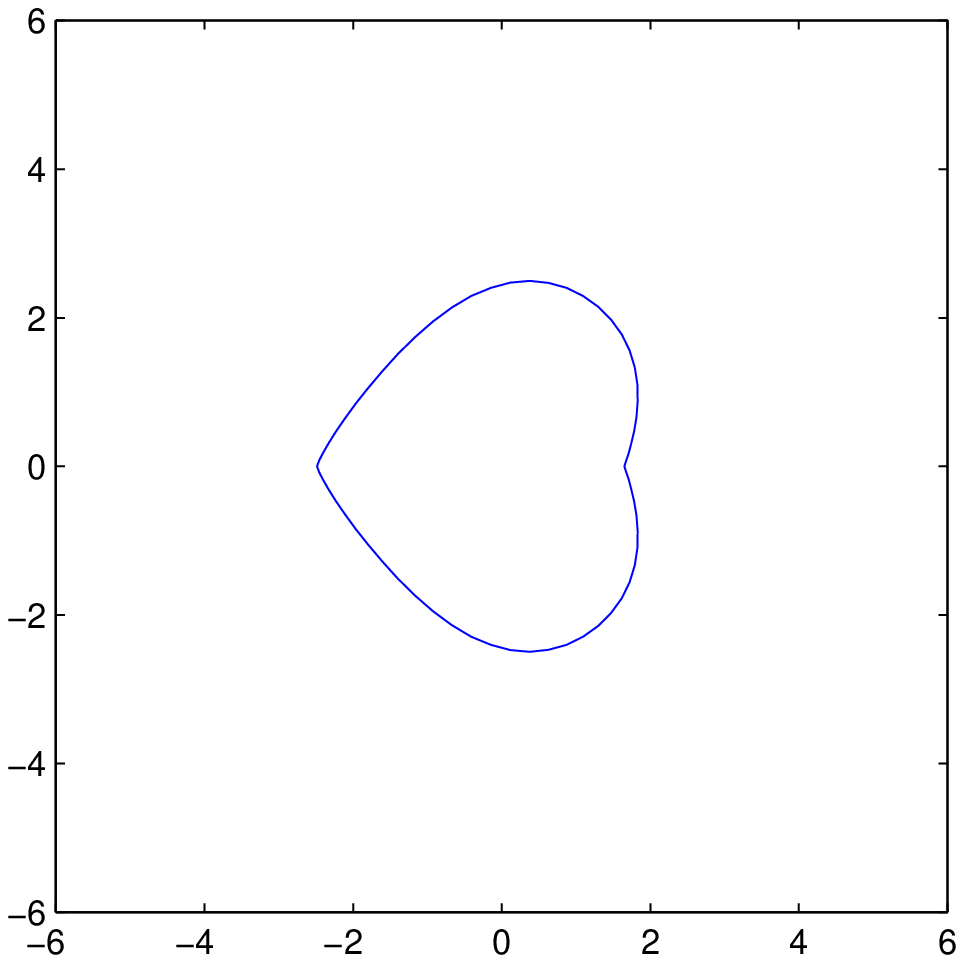}}\\
\subfigure[$t=0.6$, $H=6.6782$]{\includegraphics[width = 2in]{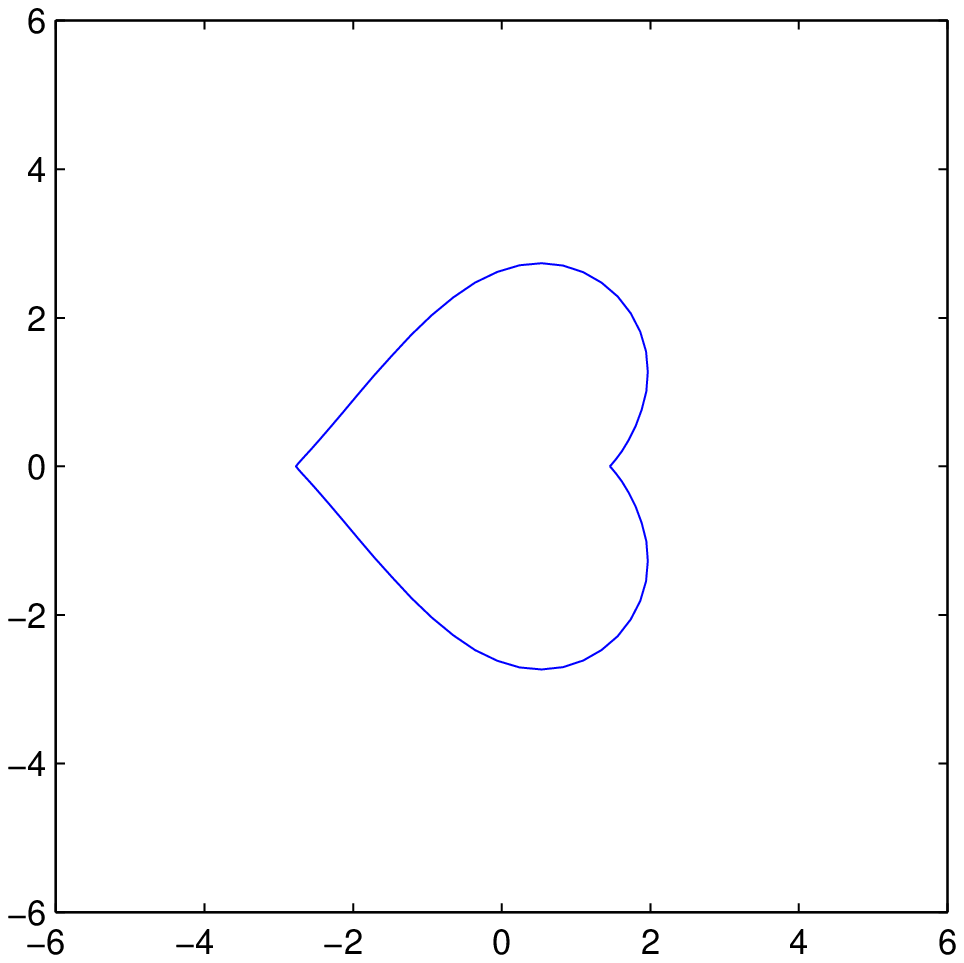}}
\subfigure[$t=0.8$, $H=11.8702$]{\includegraphics[width = 2in]{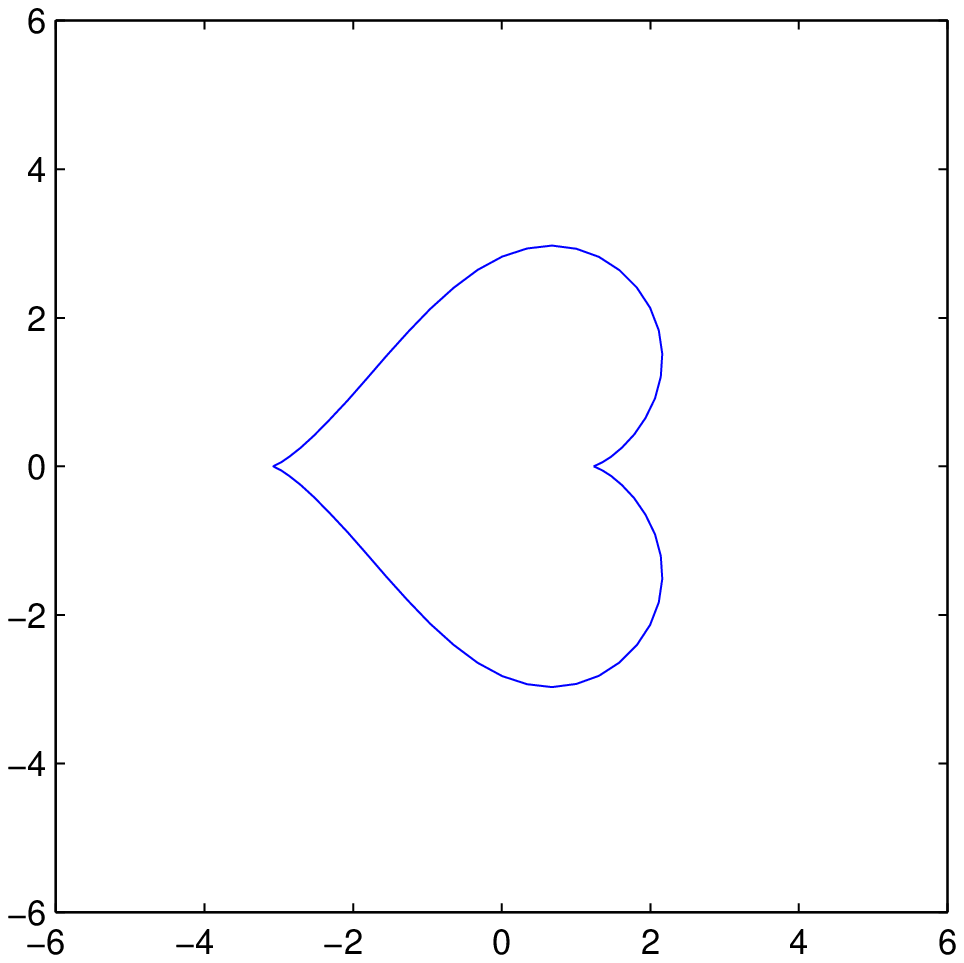}}
\subfigure[$t=1$, $H=18.5728$]{\includegraphics[width = 2in]{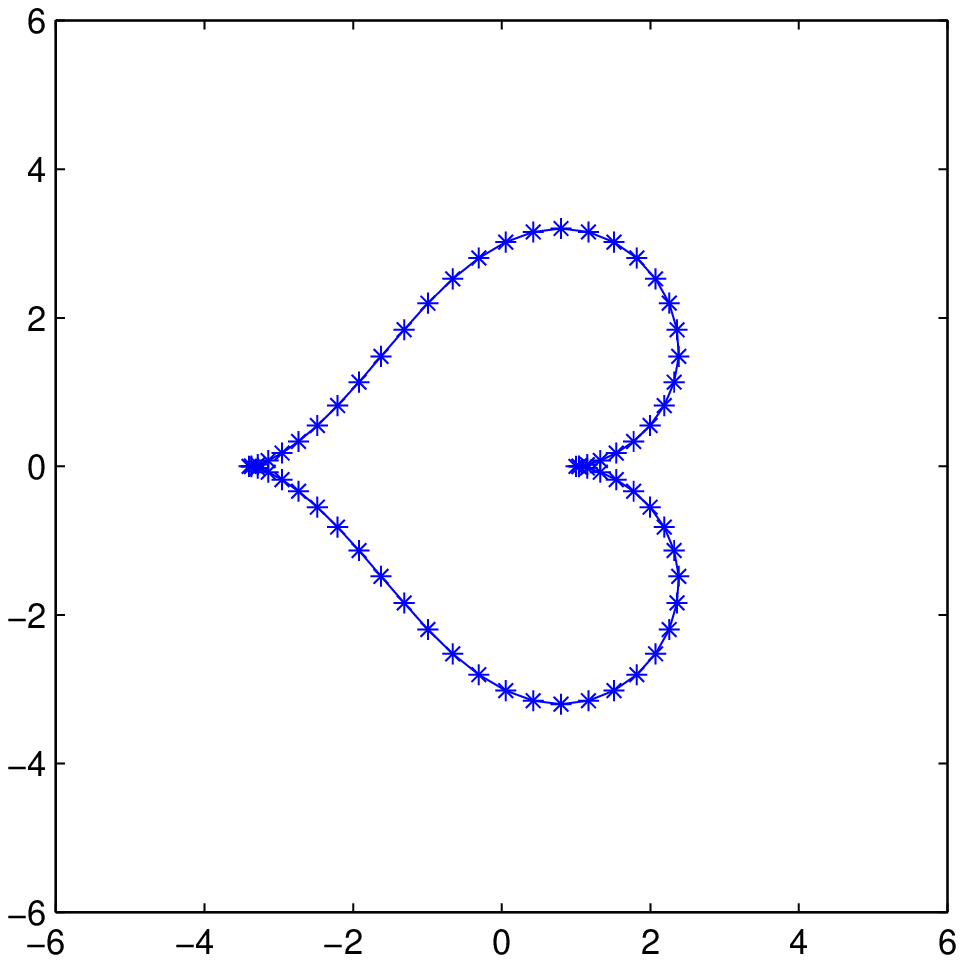}}
\caption{Diffeomorphism from a circle to a heart curve. (a) is the reference template. (b)-(e) are the interpolative shapes along the path of minimum energy found by the proposed algorithm. The dotted line in (f) is the landmark distribution of the target template, whereas the solid line is the approximate target template.}
\label{heart}
\end{figure}

\begin{table}[tbh]
\caption{The effects of different searching length $h$.}
\begin{tabular}{|l||c|c|}
\hline
\hline
$h$ & Total elapsed CPU time(sec) & \# of iteration\\
\hline
0.2 &7.5845 & 83\\ 
\hline
0.4 &  3.6137   & 40\\
\hline
0.6 & 2.3521 & 26\\
\hline
0.8 & 1.6251 & 18\\
\hline
1.0 & 3.4431 & 38\\
\hline
\hline
$>$1 & NaN &  NaN\\
\hline
\end{tabular}
\label{tab:diff-h}
\end{table}

In real applications, the choice of the search length $h$ and the number of landmarks $N$ can be subtle  for the performance of {\bf Algorithm 1}.  To investigate optimal search length, we fix the number of landmarks $N=64$ in {\bf Example 2} and evaluate the required elapsed CPU times for various $h$, given a stopping criterion. Table \ref{tab:diff-h}
suggests that the algorithm converges faster for large $h$ before the threshold, for which the algorithm begins to diverge. 
This is the typical behavior of steepest-descent type of methods. We, however, observe from our numerical experiments that the threshold for {\bf Algorithm 1} is around $h\approx 1$.

\section{Convergence property related to the metric $\mathcal{H}^{\nu}$} \label{sec:metric}

The smooth Green's kernel associated with the metric $\mathcal{H}^{\infty}$ is commonly used for template matching algorithms in the literature. In this section, we numerically investigate the convergence property of {\bf Algorithm 1}, associated with the search length $h$, the parameter $\alpha$, and the metric $\mathcal{H}^{\nu}$. In particular, our study focuses on the comparison between the smooth $\mathcal{H}^{\infty}$ and the conical $\mathcal{H}^{3/2}$ metrics. The test problem we use is again the diffeomorphism in {\bf Example 2}. Table \ref{tab:metrics} shows the iteration numbers for using various search lengths $h$ and various parameters $\alpha$ (in the Yukawa operator). Especially, (a) and (b) are the results by using the smooth Green's kernel of the metric $\mathcal{H}^{\infty}$, with the number of landmarks $N=16$ and $32$, respectively, while (c) and (d) are the same calculations with the conical Green's kernel of the metric $\mathcal{H}^{3/2}$. The ``$\times$'' mark means either the iteration number exceeds 500 before the error reaches the stopping criterion, or the algorithm blows up before 500 iterations. From the table, we observe that using the conical metric $\mathcal{H}^{3/2}$ is advantageous for fast convergence, especially for large $N$ and/or large $h$. Moreover, if we look at Table \ref{tab:metrics}(c) and (d), column by column, we observe that for the conical Green's kernel associated with the metric $\mathcal{H}^{3/2}$, the iterations are insensitive to the choice the parameter $\alpha^2$ for a fixed search length $h$. Furthermore, the table confirms the limitation of the smooth Gaussian kernel for our algorithm applied to diffeomorphisms when a large number of landmarks is required to specify the detail of a given image or shape. 
\begin{table}[h]
\caption{Numerical convergence study. (a) $N=16,\mathcal{H}^\infty$, (b) $N=32, \mathcal{H}^\infty$, (c) $N=16, \mathcal{H}^{3/2}$, (d)$N=32, \mathcal{H}^{3/2}$. The tolerance is $\epsilon=10^{-6}$.}
\label{tab:metrics}
\begin{minipage}{0.45\textwidth}
\centering
(a)
\begin{tabular}{|l||c|c|c|c|c|}
\hline
 & $h$=0.2 & 0.4 & 0.6 & 0.8 & 1.0\\
\hline
\hline
$\alpha^2$ =0.2 & 61 & 37& 32 & 48 & $\times$ \\ 
\hline
0.4 &  53   & 25 & 17 & 15 & 18\\
\hline
0.6 & 52 & 25& 19& 14& $\times$\\
\hline
0.8 & 59& 34& 22& $\times$ & $\times$\\
\hline
1.0 & 73 & 37 & $\times$ & $\times$ & $\times$ \\
\hline
\end{tabular}
\end{minipage}
\begin{minipage}{0.45\textwidth}
\centering
(b)
\begin{tabular}{|l||c|c|c|c|c|}
\hline
  & $h$=0.2 & 0.4 & 0.6 & 0.8 & 1.0\\
\hline
\hline
$\alpha^2$ =0.2 & 227 &$\times$& $\times$ & $\times$ & $\times$ \\ 
\hline
0.4 &  $\times$   & $\times$ & $\times$ &$\times$ &$\times$\\
\hline
0.6 &$\times$ &$\times$& $\times$& $\times$& $\times$\\
\hline
0.8 & $\times$& $\times$& $\times$& $\times$ & $\times$\\
\hline
1.0 & $\times$ & $\times$ & $\times$ & $\times$ & $\times$ \\
\hline
\end{tabular}
\end{minipage}

\vskip 0.5cm

\begin{minipage}{0.45\textwidth}
\centering
(c)
\begin{tabular}{|l||c|c|c|c|c|}
\hline
  & $h=0.2$ & 0.4 & 0.6 & 0.8 & 1.0\\
\hline
\hline
$\alpha^2$ =0.2 & 39 & 18& 12 & 9 & 8 \\ 
\hline
0.4 &  40   & 18 & 11 & 8 & 9\\
\hline
0.6 & 40 & 18& 11& 9& 11\\
\hline
0.8 & 40& 18& 11& 9 & 12\\
\hline
1.0 & 40 & 18 & 11 & 10 & 13 \\
\hline
\end{tabular}
\end{minipage}
\begin{minipage}{0.45\textwidth}
\centering
(d)
\begin{tabular}{|l||c|c|c|c|c|}
\hline
  & $h$=0.2 & 0.4 & 0.6 & 0.8 & 1.0\\
\hline
\hline
$\alpha^2$ =0.2 & 46 & 22& 16 & 15 & $\times$ \\ 
\hline
0.4 &  47   & 22 & 13 & 10 & 17\\
\hline
0.6 & 54 & 25& 16& 11& 20\\
\hline
0.8 & 58& 27& 17& 12 & 21\\
\hline
1.0 & 60 & 29 & 18 & 12 & 24 \\
\hline
\end{tabular}
\end{minipage}

\end{table}

\begin{figure}[tbh]
\subfigure[$N=16,\mathcal{H}^{\infty}$]{\includegraphics[width = 3.in]{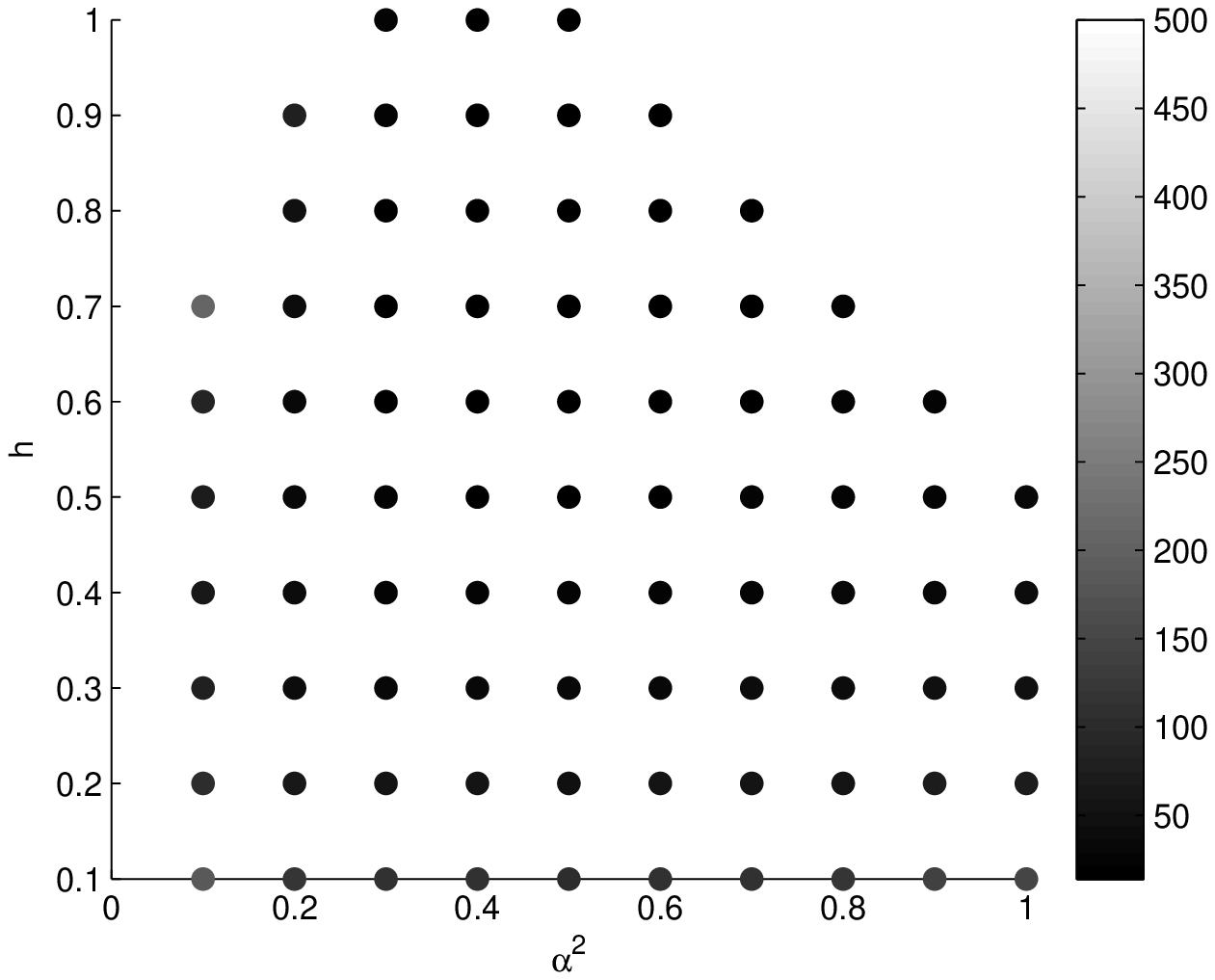}} 
\subfigure[$N=32,\mathcal{H}^{\infty}$]{\includegraphics[width = 3.in]{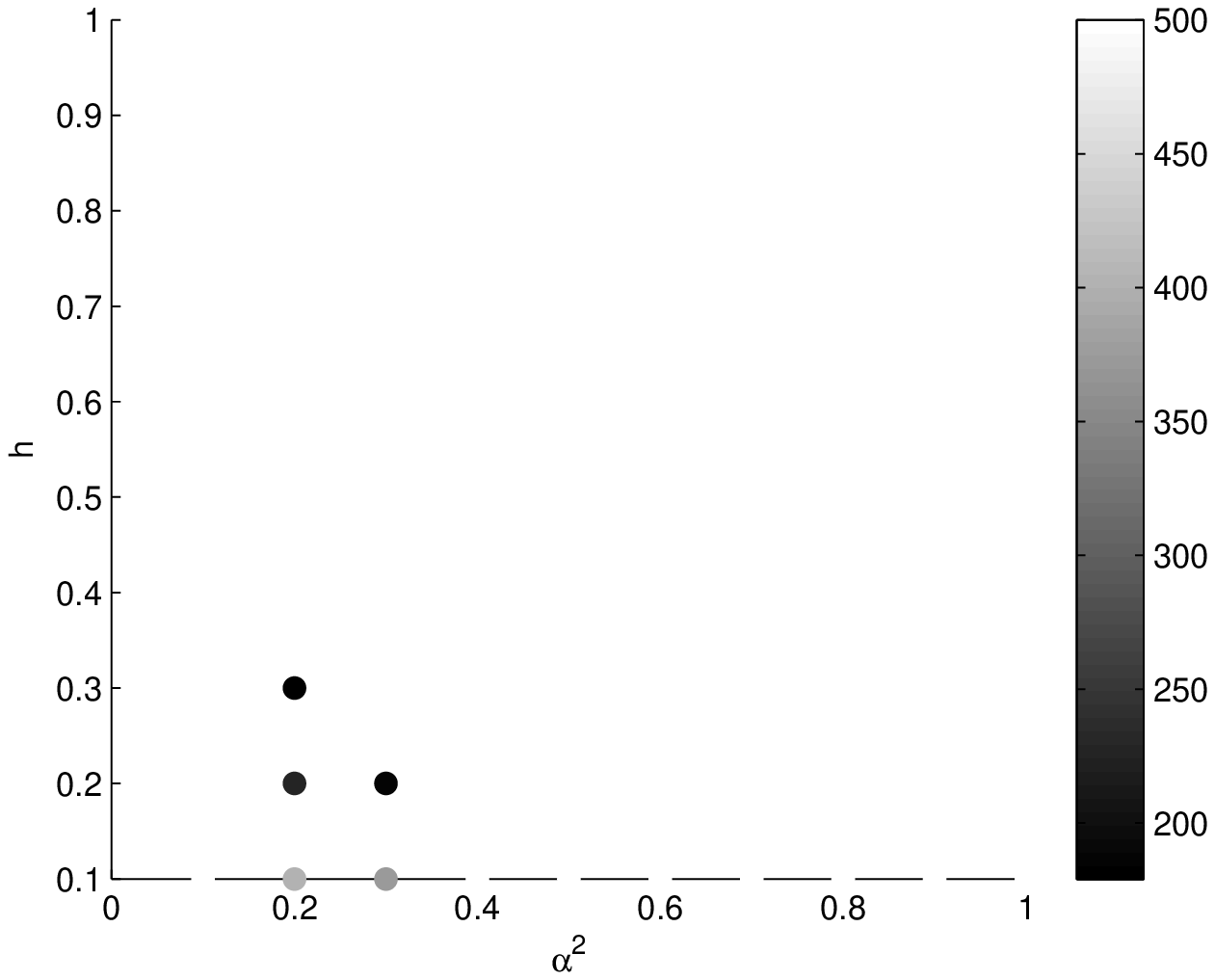}}\\
\subfigure[$N=16,\mathcal{H}^{3/2}$]{\includegraphics[width = 3.in]{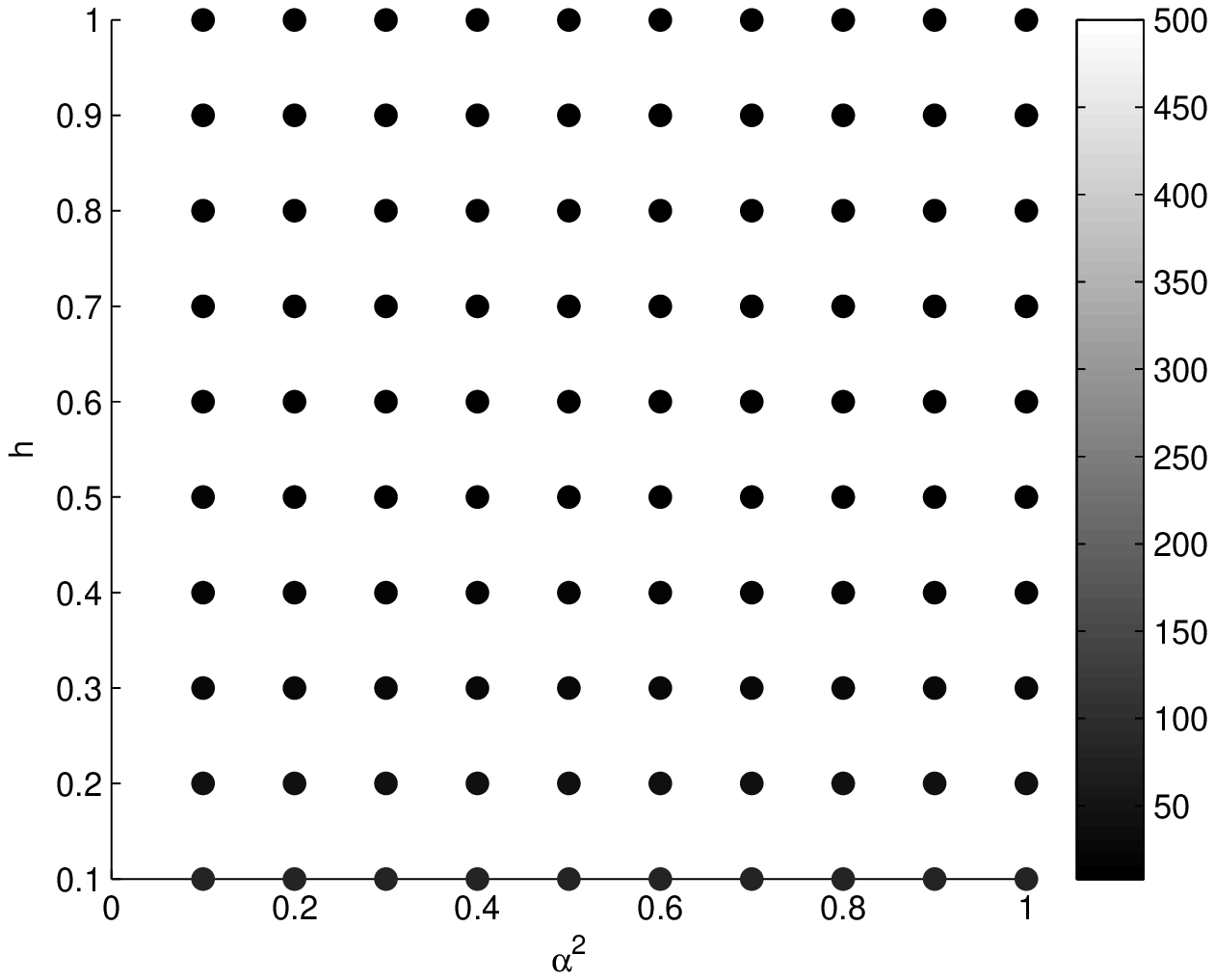}} 
\subfigure[$N=32,\mathcal{H}^{3/2}$]{\includegraphics[width = 3.in]{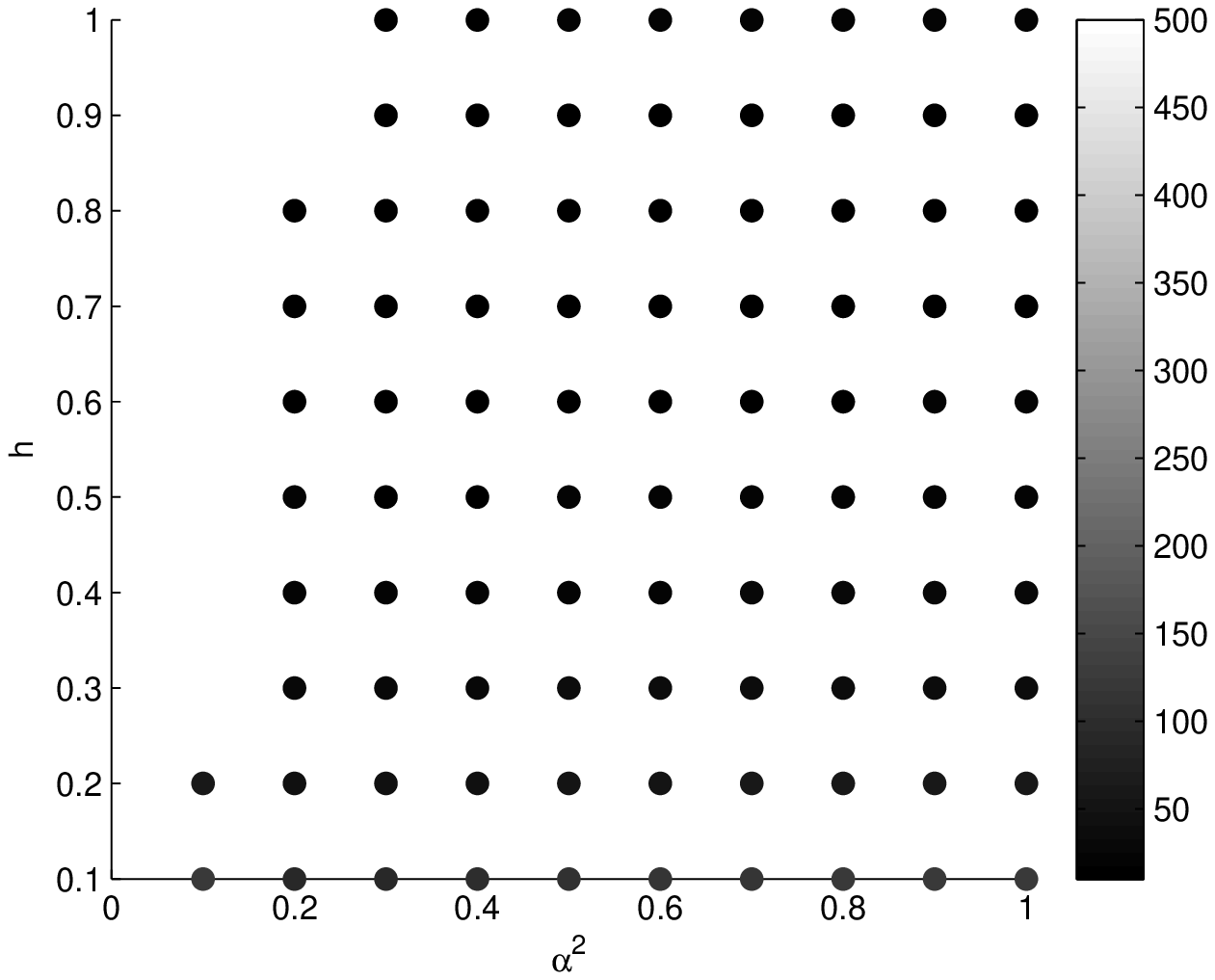}}
\caption{Numerical convergence study of the metrics $\mathcal{H}^{\infty}$ and $\mathcal{H}^{3/2}$ in the $\alpha^2\times h$ parameter space. The tolerance is $\epsilon=10^{-6}$ in the maximum norm and the gray level represents the number of iterations.}
\label{fig:metric_comp}
\end{figure}

To better visualize the comparison, we place the convergence results in the $\alpha^2\times h$ parameter space, using a grid size $0.1 \times 0.1$. A black dot  indicates convergence, while a white dot suggests no convergence. The darker the dot is, the less iterations it takes to reach the tolerance. Figure \ref{fig:metric_comp} supports our previous conjecture that particles have a shorter interaction range with the non-smooth Green's kernel. As a result, the algorithm is more stable and robust when the  non-smooth Green's kernel is used, especially when the number of landmarks (particles) is large.

\section{The prediction ability of {\bf Algorithm 1}}
\begin{figure}[tbh]
\includegraphics[width=3.7in]{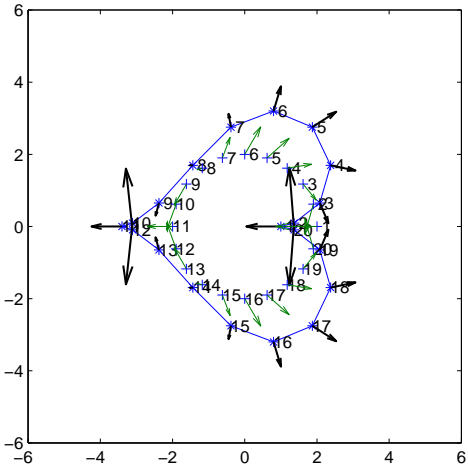}
\caption{The approximate momenta in {\bf Example 2}.}
\label{fig:momenta}
\end{figure}

One of the advantages of  {\bf Algorithm 1} is its ability of predicting the deformation beyond the target template.  Traditionally, algorithms for template matching solve the Euler-Lagrangian equation of the optimization problem between two templates and obtain the deformations of the reference template through the iterative process. {\bf Algorithm 1}, however, approximates the initial momenta of the particle system that drive the deformation from template $A$ to template $B$. As a result, the algorithm could continue the deformation of the reference template beyond the target template by using the approximate momenta. Figure \ref{fig:momenta} shows the approximate initial momenta (thin arrows) that drive the diffeomorphism from the circle to the $4^{th}$-order hear curve, and the momentum of each landmark at the target contour (thick arrows) calculated from the approximate initial momenta. The continuation of the $4^{th}$-order heart curve, based on the found approximate initial momenta, is easily predicted. 

\begin{figure}[tbh]
\includegraphics[width=5.4in]{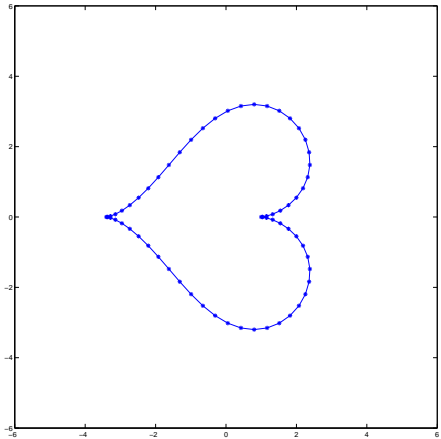}
\caption{Comparison between the prediction (dots) and the actual data (solid line).}
\label{fig:predict}
\end{figure}

We further investigate the prediction ability of {\bf Algorithm 1} by the following experiment.
\begin{enumerate}
\item Take phase Figure \ref{heart}(d) as our target template, denoted as $I_{0.6}$;
\item Use {\bf Algorithm 1} to  approximate the initial momenta (or the flow $\bs{u}$) that drive $I_0$ to $I_{0.6}$, i.e. $I_{0.6}={\phi}(\bs{u};I_0,0.6)$;
\item Use the $\bs{u}$ obtained from step (2) to evolve the particle system upto $t=1$ for $I_0$, i.e $\hat{I_1}={\phi}(\bs{u};I_0,1)$;
\item Compare  $\hat{I_1}$ to the original target template $I_1$ given by Eq. (\ref{eq:heart}) (Figure \ref{heart}(f)).
\end{enumerate}
Figure \ref{fig:predict} plots $\hat{I_1}$ and $I_1$ on top of each other. The solid line is Eq. (\ref{eq:heart}) and the dots are landmarks computed by the above steps (1) -(3), using {\bf Algorithm 1}. We see that they are visually indistinguishable. This consistency experiment shows that {\bf Algorithm 1} correctly predicts the deformation beyond the target template.

\section{Applications of shape analysis}\label{sec:app}
{\bf Algorihtm 1} can be applied to shape and cluster analysis.  For this type of problems, a certain measure needs to be defined and computed among a given set of patterns against multiple references for the purpose of classification. The Hamiltonian $H$ is an example for such a measure.


Before introducing {\bf Algorithm 1} and the Hamiltonian $H$  as tools for shape and cluster analysis, we like to point out some properties about $H$. One is that $H $ is a semi-metric defined in the tangent space of a manifold other than in the manifold itself. Points that are close to each other in a manifold can be pulled apart in the tangent space if the curvature is positive, while some far distant sets of points may appear as a cluster in the tangent space \cite{bib:mumford}. 
Suppose $A$ and $B$ are two targets templates, and $C$ is the reference template, if the difference $|H(A,C)-H(B,C)|$ is small, it is not necessary to imply that $H(A,B)$ is small. Here $H(A, C)$ represents the Hamiltonian between $A$ and $C$. Geometrically, since the Hamiltonian is an indication of how easy (or difficult) it is to deform from one planar curve to the other,  a small value of this measurement does not necessarily imply that the two curves is geometrically similar. This is somewhat against our intuition, as shown in our next example.

\begin{figure}[tbh]
\subfigure[$H=0$ ]{\includegraphics[width = 2in]{fig1.eps}} 
\subfigure[$H=7.2456$]{\includegraphics[width = 2in]{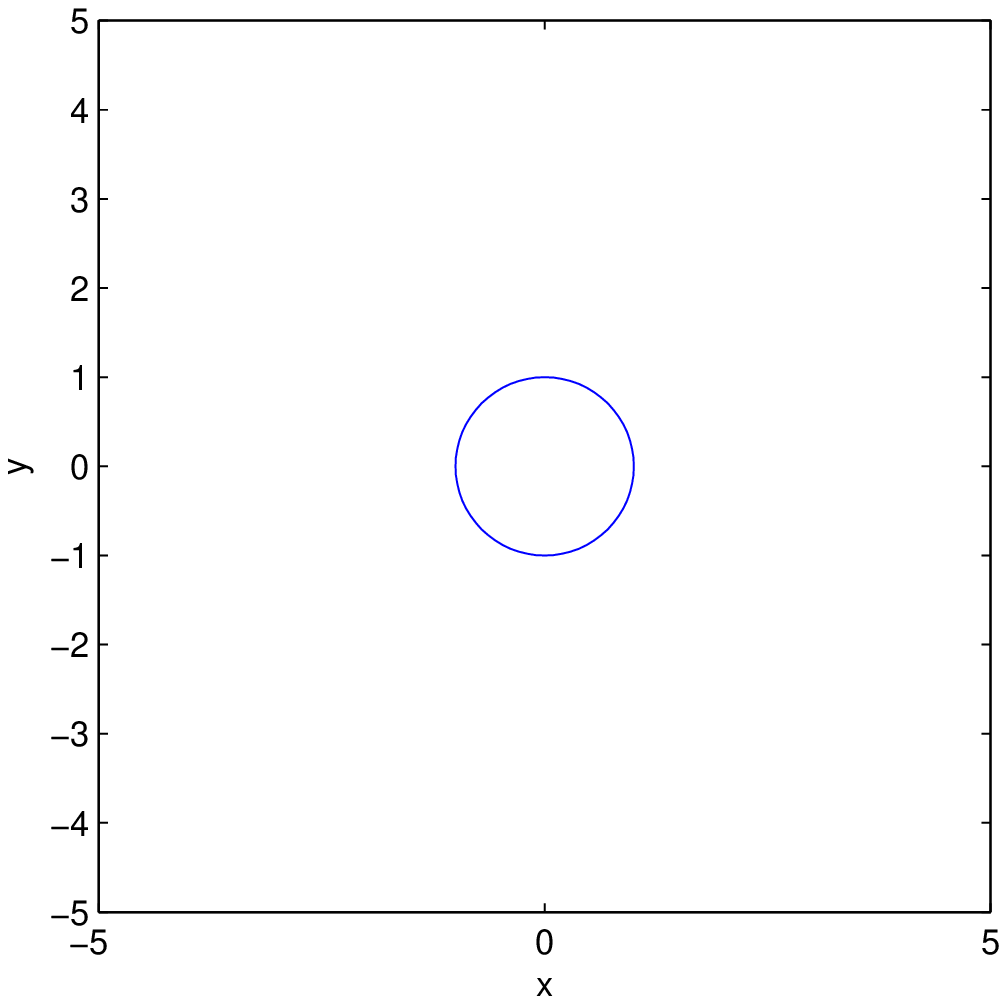}} 
\subfigure[$H=9.1209$]{\includegraphics[width = 2in]{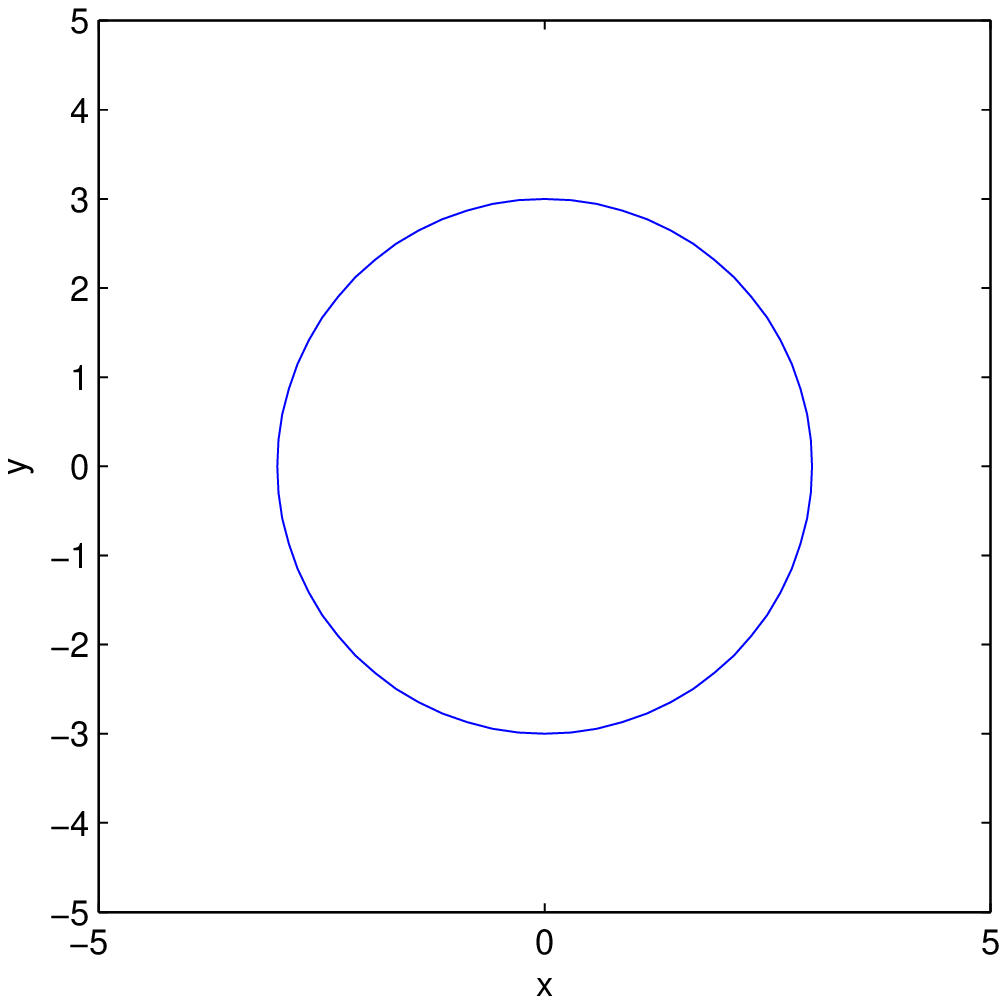}}\\
\subfigure[$H=19.2438$]{\includegraphics[width = 2in]{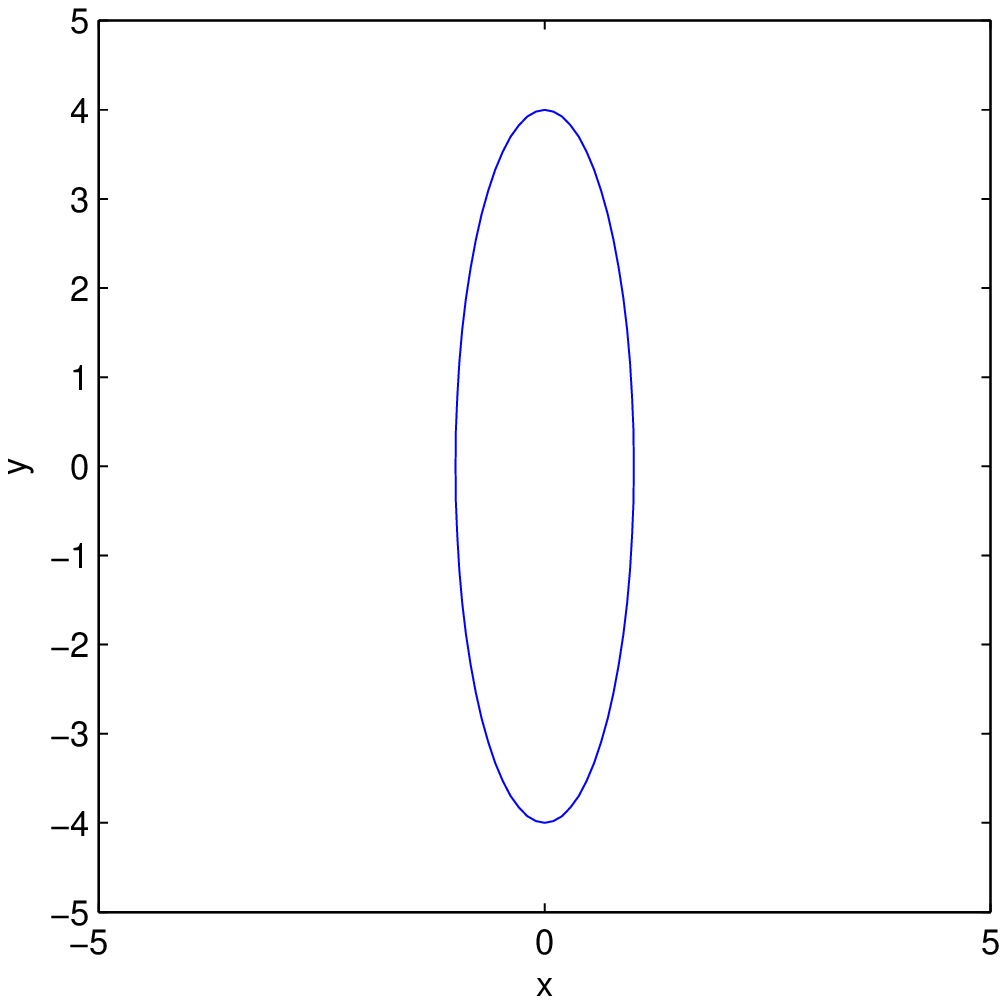}} 
\subfigure[$H=8.2542$]{\includegraphics[width = 2in]{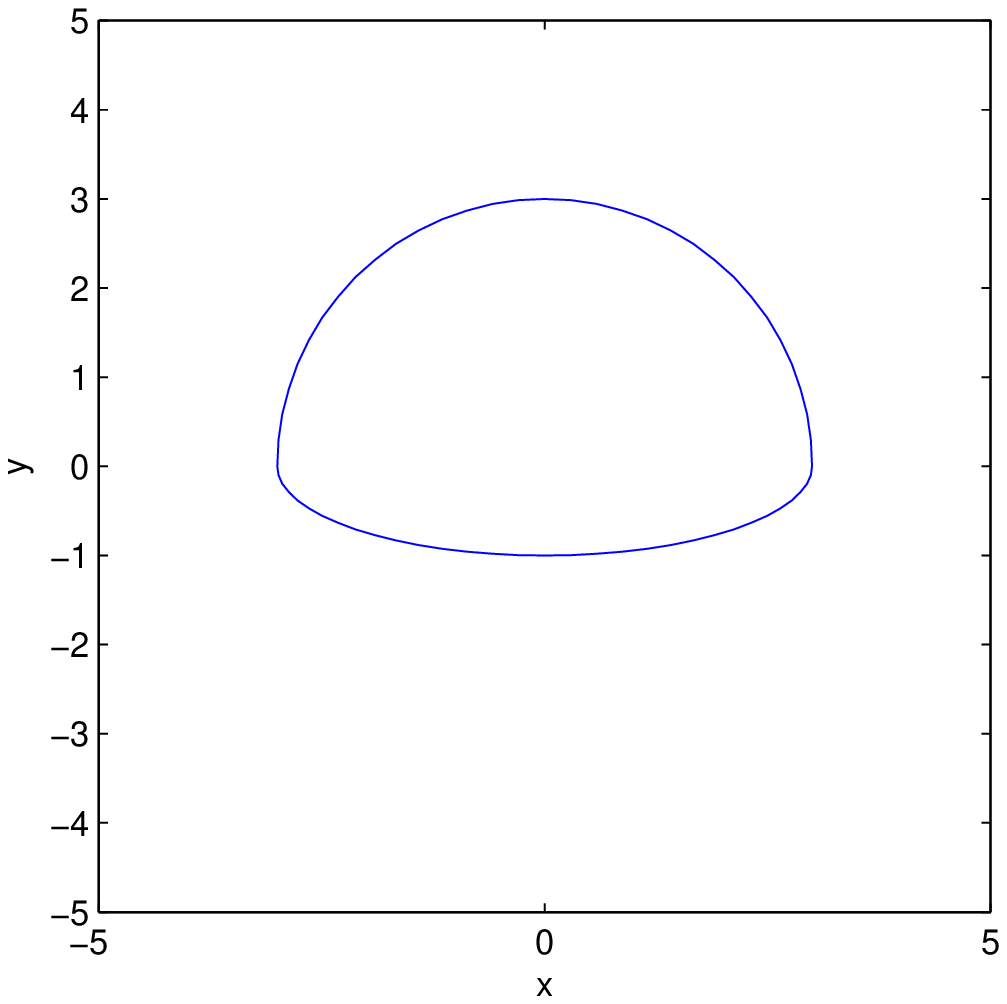}} 
\subfigure[$H=2.9958$]{\includegraphics[width = 2in]{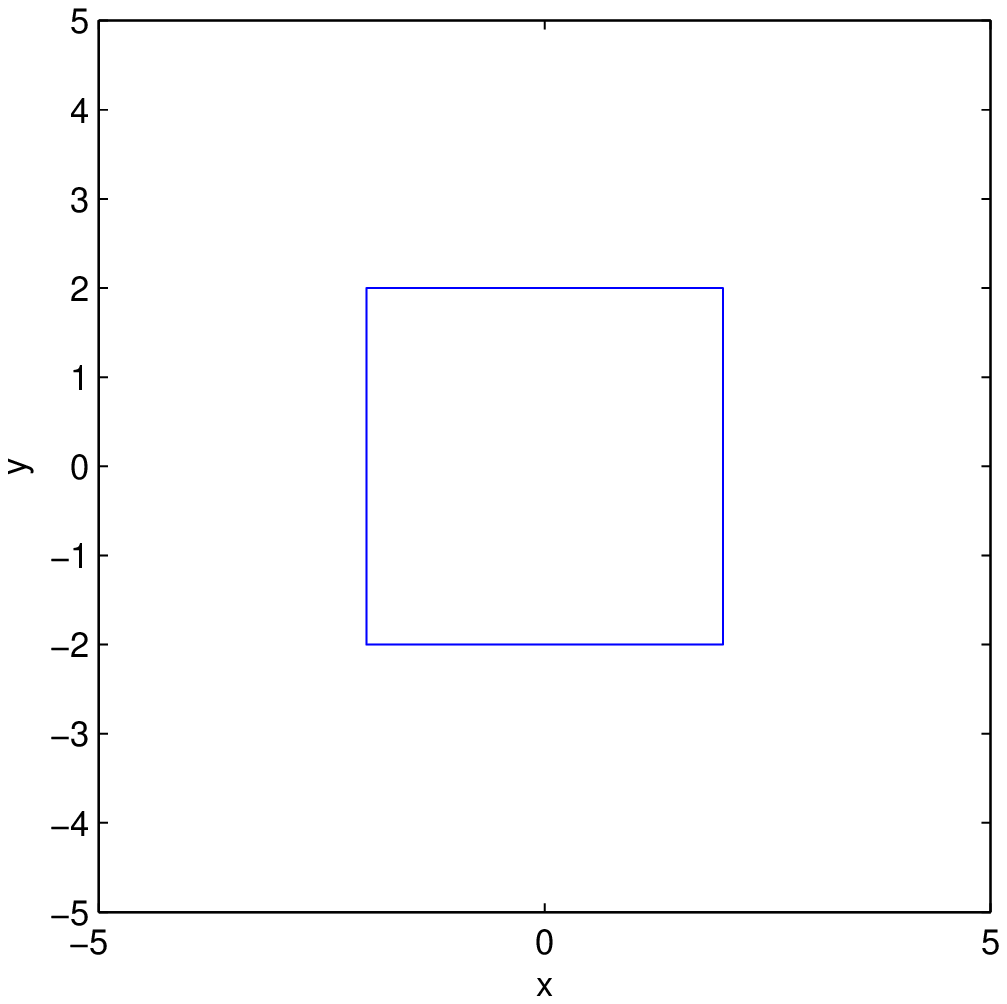}}\\
\caption{Examples of the Hamiltonian.} 
\label{fig:H_distance}
\end{figure}

In Figure \ref{fig:H_distance}, (a) is a circle with $r=2$, (b) is a shrunken circle with $r=1$, (c) is the dilated circle with $r=3$, (d) is an ellipse with $a=1, b=4$, (e) is the closed curve with upper half from a circle($r=3$) and lower half from an ellipse$(a=3,b=1)$, and (f) is a square with side length $l=4$. Suppose that (a) is our reference template. The Hamiltonian from reference (a) to targets (b)-(f), in the ascent order, is $H(af)< H(ab)< H(ae)< H(ac)<H(ad)$. $N=64$ for all calculations. Our intuition about a metric distance is that a shorter metric distance implies a closer geometric similarity. If this ought to be true, then the square (f) is ``closest`` to the circle (a) among the five patterns.  This is obviously against our intuition. Nevertheless, the results could be understood by the following properties of the Hamiltonian.

First of all, $H(A,B)\neq H(A,\mathcal{G}(B))$ for $\mathcal{G}$ in group generated by dilation (or contraction), rotation and translation. This suggests that any dilations (or contractions), rotations or translations will increase the Hamiltonian. In Figure \ref{fig:H_distance}, (b) is a contraction of (a) while (c) is a dilation of (a), which results in large Hamiltonian for (a)\&(b) and (a)\&(c). Also, in Figure \ref{fig:H_distance}, except (e) the centroids of the curves are located at the origin, which causes a large Hamiltonian between (a) and (e).
We observe that in Figure \ref{fig:H_distance}, the volume enclosed by the contour curves plays an important role to determine the Hamiltonian. This explains why the Hamiltonian from (f) to (a) is small, since the enclosed area of (f) is close to that of (a). The Hamiltonian carries other information, too. For example, similar to (a), the contour (d) also has an enclosed area that is close to (a), but (d) has the furthest distance, whereas (f) has the closest one. This can be explained by the symmetry around the centroid (geometric center).  Because contours (a) and (f) are both symmetric around the centroid, the distance between (a) and (f) is shorter than the distance between (a) and (d) (contour (d) is only symmetric about $y$-axis). It is worth pointing out that the Hamiltonian has a nice property: $H(A,B) = H(\mathcal{G}(A),\mathcal{G}(B))$ for $\mathcal{G} \in \text{ISO}(\mathbb{R}^n)$, where $\text{ISO}(\mathbb{R}^n)$ denotes the linear isometric transformation group of $\mathbb{R}^n$, i.e. rotation, reflection and translation. TABLE \ref{tab:isoH} shows a numerical example to illustrate this property. This property brings up the potential of the Hamiltonian as a tool for shape analysis.

\begin{table}[tbh]
\caption{An example: ISO invariance of H. Here A is shape (a), B is shape (e), $\Delta H$ calculates $H(\mathcal{G}(A),\mathcal{G}(B))-H(A,B)$}
\label{tab:isoH}
\begin{tabular}{cc}
\hline
Transformation $\mathcal{G}$ & \qquad Change of $H(\Delta H)$\\
\hline
counterclockwise rotation $\pi/3$ &  \qquad -2.3e-14\\
\hline
translation along (0.5, 0.5)  &   \qquad -7.1e-15\\
\hline 
reflection with respect to $\hat{\bs{x}}$  &\qquad 0\\
\hline
\end{tabular}
\end{table}

Although the Hamiltonian only indicates how much energy is consumed in a diffeomorphism, and the geometric similarity between two templates does not always lead to a small distance between them, vise versa, we will show in the following subsection that {\bf Algorithm 1} and the Hamiltonian $H$ compose a tool for image registration or shape analysis. Finally, we remark that the Hamiltonian is symmetry between the reference and the target templates. For example, in Figure \ref{fig:H_distance} our calculations show that the Hamiltonian $H(ec)=10.0070$ ((e) is the reference and (c) is the target), whereas $H(ce)=10.0132$ ((c) is the reference and (e) is the target).


\subsection{Example of using {\bf Algorithm 1} for shape analysis or patter recognition}  We demonstrate an example of using {\bf Algorithm 1} and the Hamiltonian $H$ for shape analysis or patter recognition. Suppose that we want to determine whether two shapes $A$ and $B$ belong to the same cluster. We can perform the following algorithm for the determination.

\begin{enumerate}
\item Select a reference template $C$. 
\item Apply {\bf Algorithm 1} for template matching between $C$ \& $A$ and $C$ \& $B$.
\item When the matching is complete, compute the Hamiltonian $H(CA)$ and $H(CB)$.
\item Select a second template $D$, where the Hamiltonian $H(CD)$ is large.
\item Repeat steps (2) and (3) for $H(DA)$ and $H(DB)$. 
\item Apply {\bf Algorithm 1} for template matching between $A$ \& $B$ and compute the Hamiltonian $H(AB)$.
\item If $| H(CA)- H(CB) |$, $| H(DA)- H(DB) |$ and $H(AB)$ are all smaller than the chosen thresholds, then it is possible that $A$ and $B$  belong to the same cluster, otherwise they may not.
\end{enumerate} 
\begin{figure}[tbh]
\subfigure[Target $A$]{\includegraphics[width = 2.5in]{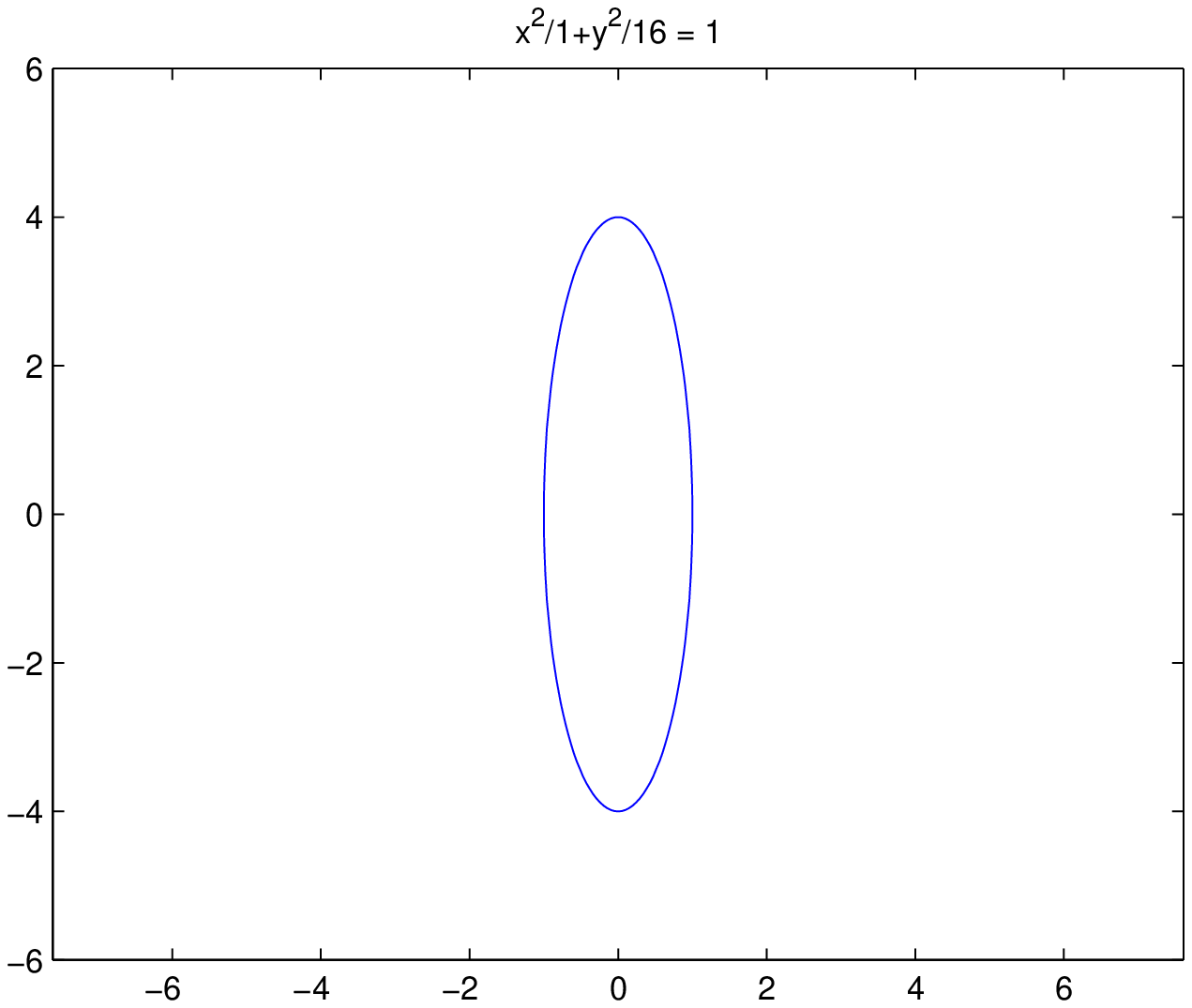}} 
\subfigure[Target $B$]{\includegraphics[width = 2.5in]{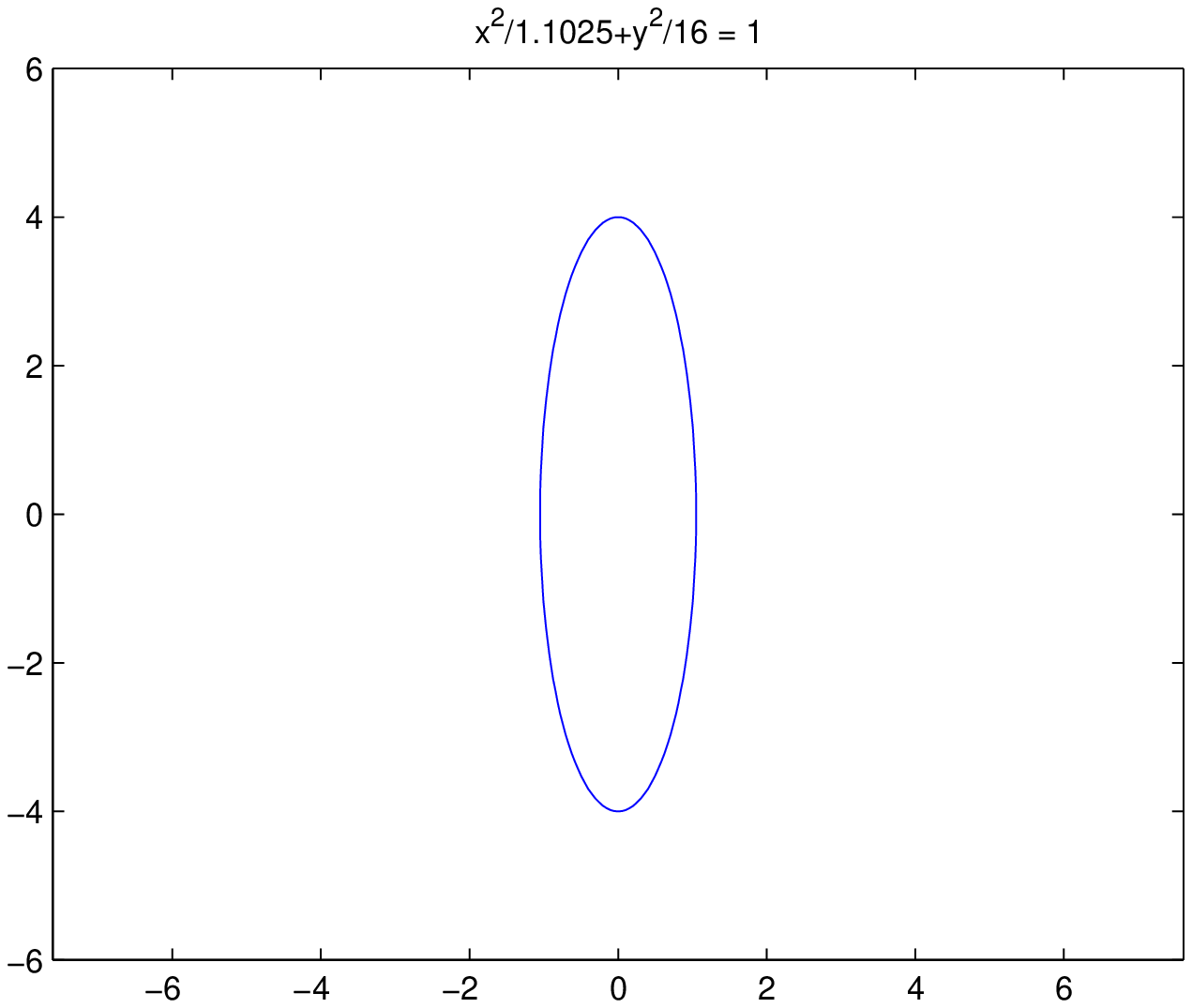}}\\
\subfigure[Reference $C$]{\includegraphics[width = 2.5in]{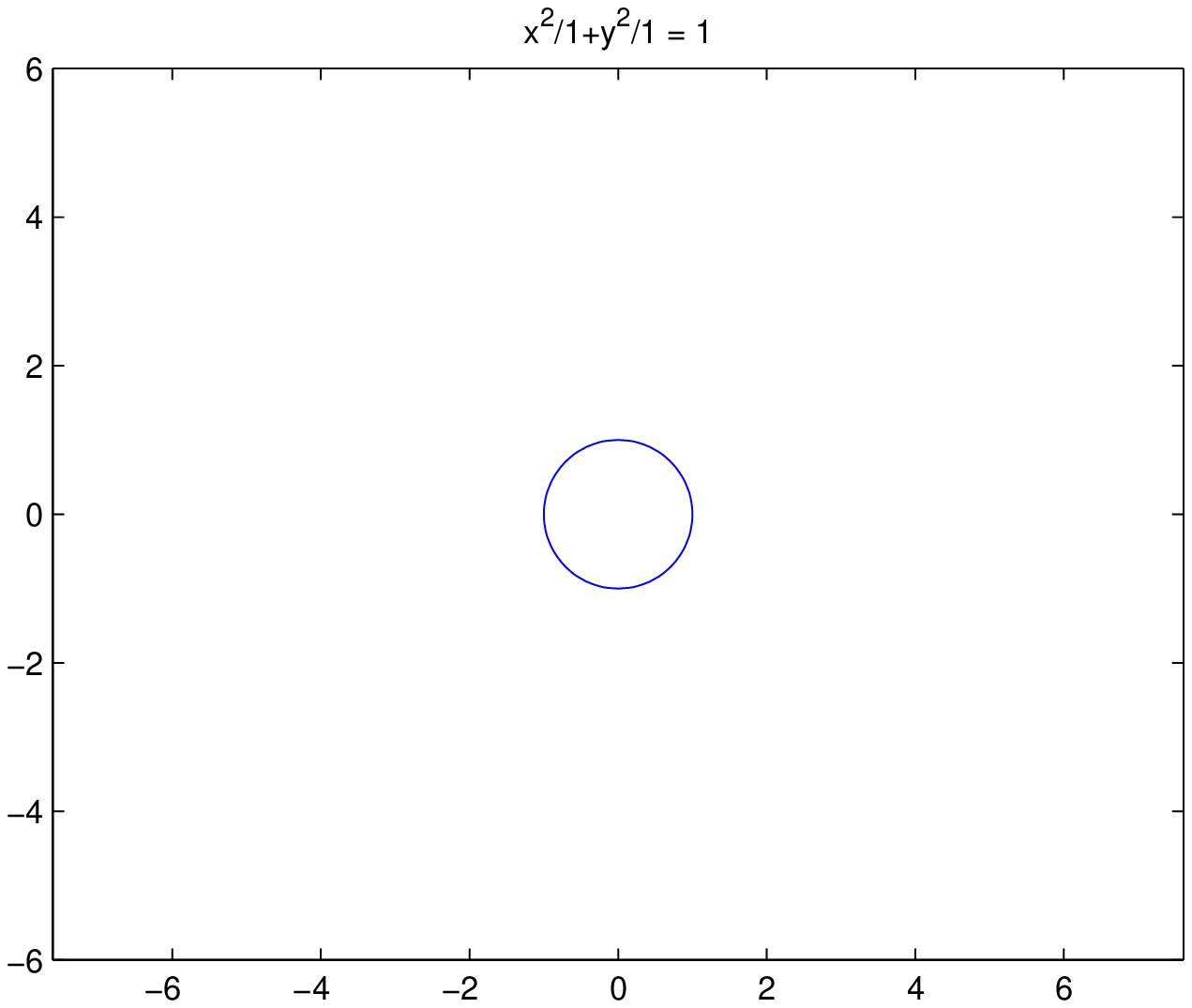}} 
\subfigure[Reference $D$]{\includegraphics[width = 2.5in]{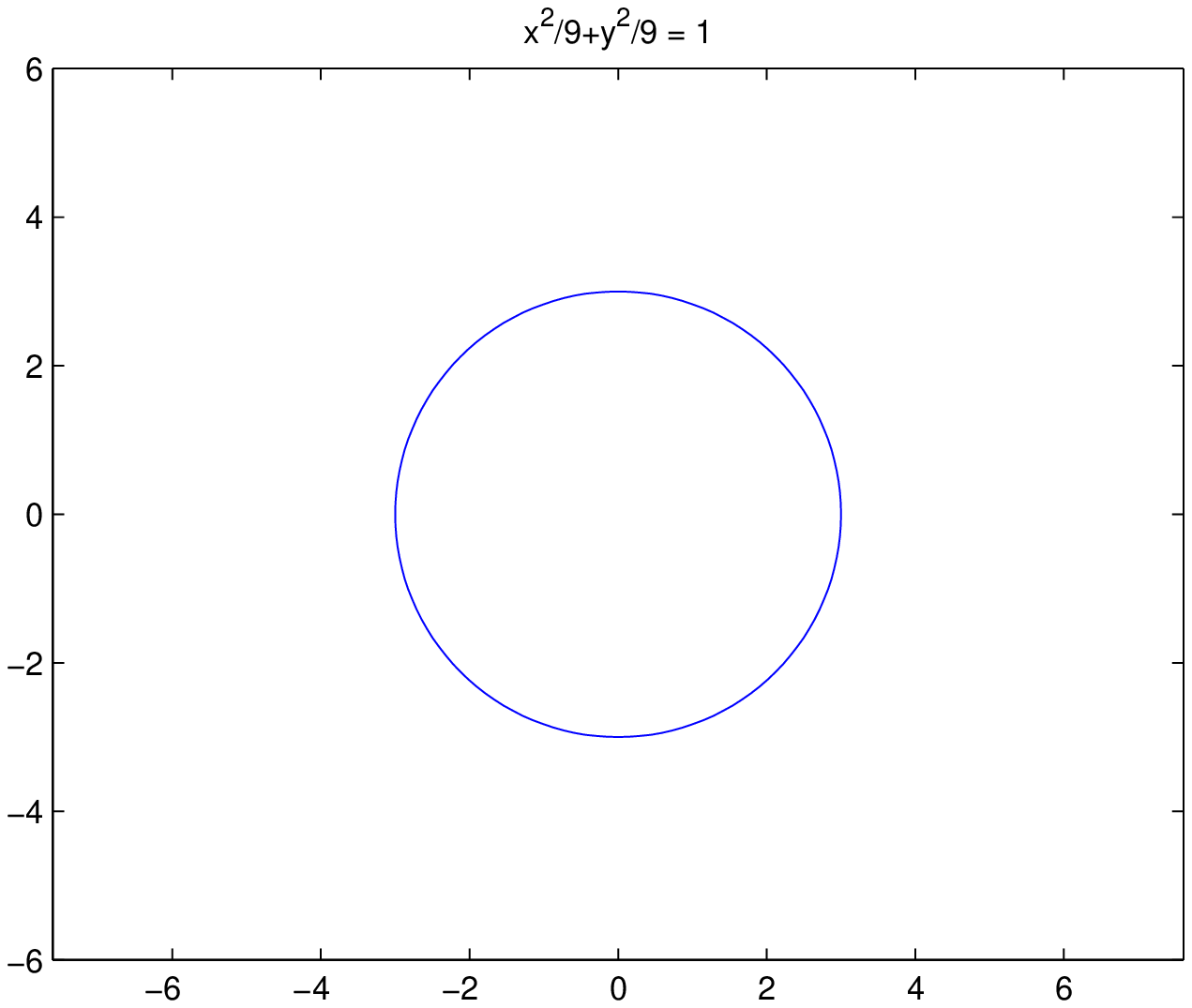}}
\caption{Targets and references for the example of using {\bf Algorithm 1} and the Hamiltonian $H$ for the shape analysis.}
\label{fig:shapes}
\end{figure}
In Figure \ref{fig:shapes}, the target shapes $A$ and $B$ are ellipses with major axes both equal to 4. The minor axes are 1 and 1.05, respectively. The goal is to determine whether $A$ and $B$ belong to the same shape cluster.

We choose the references $C$ and $D$ that are circles with radii equal to 1 and 3, respectively. The Hamiltonian between $C$ and $D$ is $H(CD)=9.1209$. Following the above algorithm, we apply {\bf Algorithm 1} for template matching between the references and the targets. We compute the following Hamiltonian: $H(CA)=19.2438$, $H(CB)=18.7700$, $H(DA)=28.3570$, $H(DB)=27.0526$, and $H(AB)=0.0215$. Therefore, we obtain $|H(CA)-H(CB)|=0.4743$, $|H(DA)-H(DB)|=1.3044$. If we require $H(AB)\le 0.05$, $|H(CA)-H(CB)|\le 1.5$, and $|H(DA)-H(DB)|\le 1.5$ for $A$ and $B$ being in the same cluster, then we can classify these two ellipses into the same cluster. On the other hand, if for being in the same cluster means $|H(CA)-H(CB)|\le 1$ and $|H(DA)-H(DB)|\le 1$, then $A$ and $B$ are not in the same cluster. 

The above example involves some very simple geometry, and in fact it is easy to visually make the classification, if we want. In practice, however, shapes in a library are usually irregular, complex, and in large quantity. It makes sense to use an efficient algorithm for clustering analysis. We remark that it may be necessary to preprocess the targets to eliminate the effects of relative translaitons, rotations, and dilations (contractions) for a better classification \cite{bib:Bookstein}.  We also remark that when using {\bf Algorithm 1} and the Hamiltonian for shape analysis or pattern recognition, it is always a good strategy to use multiple reference templates ($\ge 2$)  before making the classification. 
\begin{note}
In this numerical example, we found that the Hamiltonian numerically satisfies the triangle inequality for all cases, for example, $H(A, B) \le H(A, C)+ H(C, B)$. This is an indication that $H$ could be a metric, but we do not pursue the proof in this paper.
\end{note}

\section{Numerical Experiments - Inexact Matching}
In this section, we repeat the previous examples, {\bf Example 1} and  {\bf Example 2}, in Section \ref{sec:exact}, but place them in the context of inexact matching. As pointed out in \cite{bib:Holm, bib:eqs}, $\sigma^2$ in the particle system for inexact matching (see Eq. (\ref{eq:N-particle-inexact})) provides a slightly inexact advection to move ahead the deformation and avoid the formation of singularities \cite{bib:hrty04}. This kind of singularities can occur for smooth Gaussian kernels, since smooth Gaussians possess capture orbits when two particles are having in-line collision \cite{NP, bib:ckl14}. The non-smooth Green's kernel ($\nu=3/2$) used in this paper, on the other hand, does not have capture orbits \cite{bib:ckl14}. Its scattering orbits allow particles to move pass each other (crossing).  Although in evolution of biological structures, crossing is rare and we usually do not expect this to happen, the non-smooth Green's kernel, nonetheless, has the advantage of simulating this kind of natural process under the setting of exact matching.

The numerical experiments of inexact matching in this section, however, provide a different perspective of matching, other than preventing the formation of singularities. Suppose that the target landmarks used for the matching are approximates to the exact locations.  For this perspective, by choosing the size of $\sigma^2$, we allow the reference converge to something that is different from the target.   The only modification required in {\bf Algorithm 1}, other than using the inexact particle system, Eq. (\ref{eq:N-particle-inexact}), for this perspective is that the stopping criterion is chosen to be a small  difference of the approximate initial momenta computed between the $k^{th}$ and $(k -1)^{th}$ iterations.

\begin{figure}[tbh]
\subfigure[$t=0$]{\includegraphics[width = 2in]{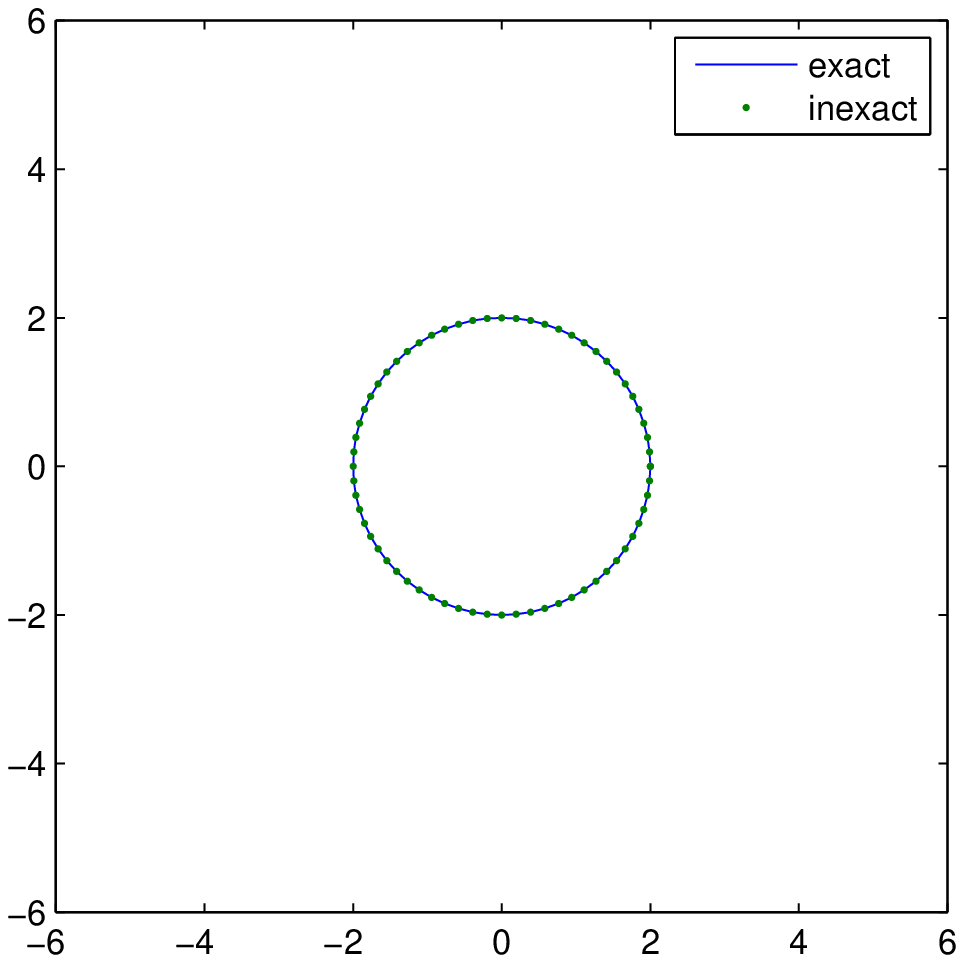}}  
\subfigure[$t=0.2$]{\includegraphics[width = 2in]{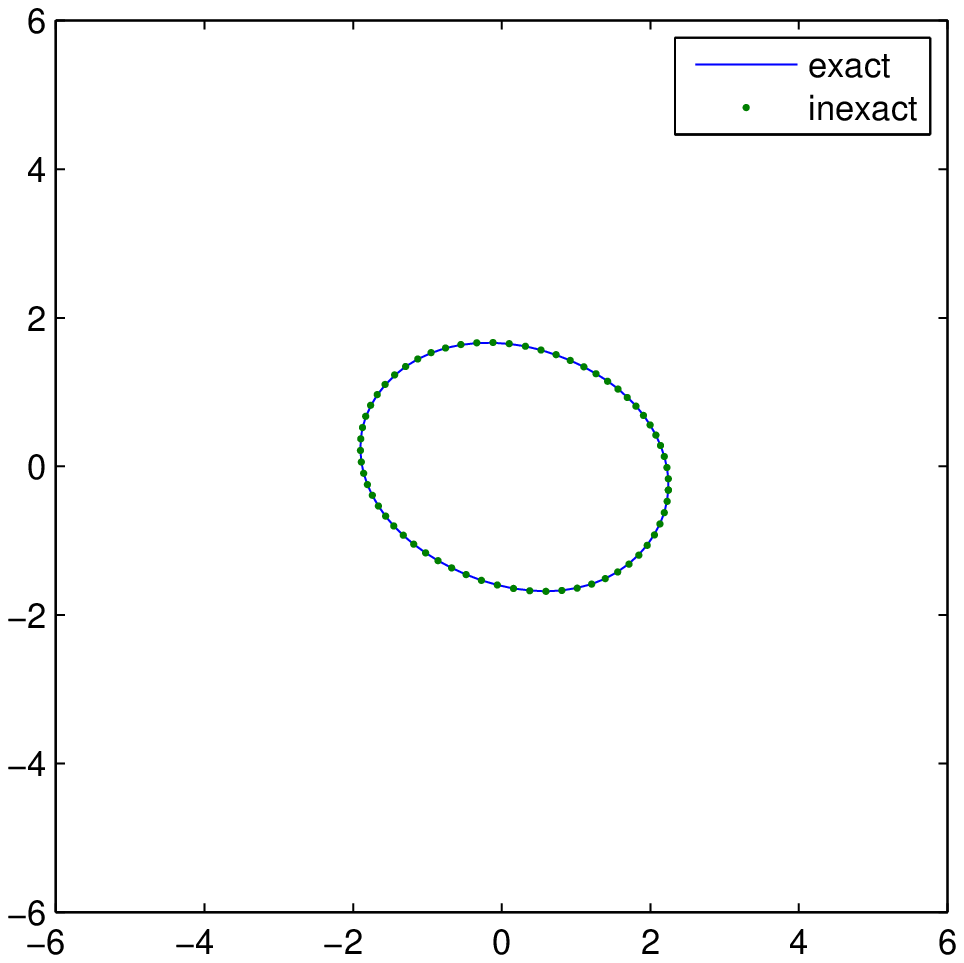}}
\subfigure[$t=0.4$]{\includegraphics[width = 2in]{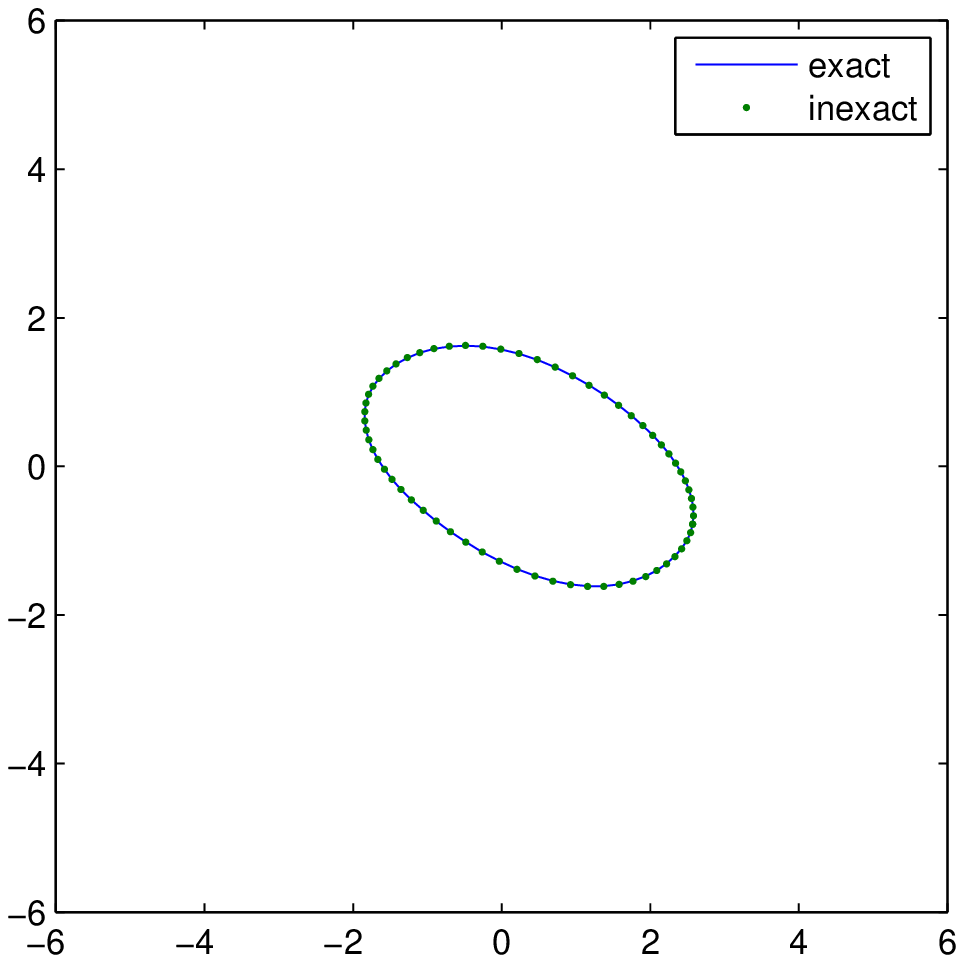}}\\  
\subfigure[$t=0.6$]{\includegraphics[width = 2in]{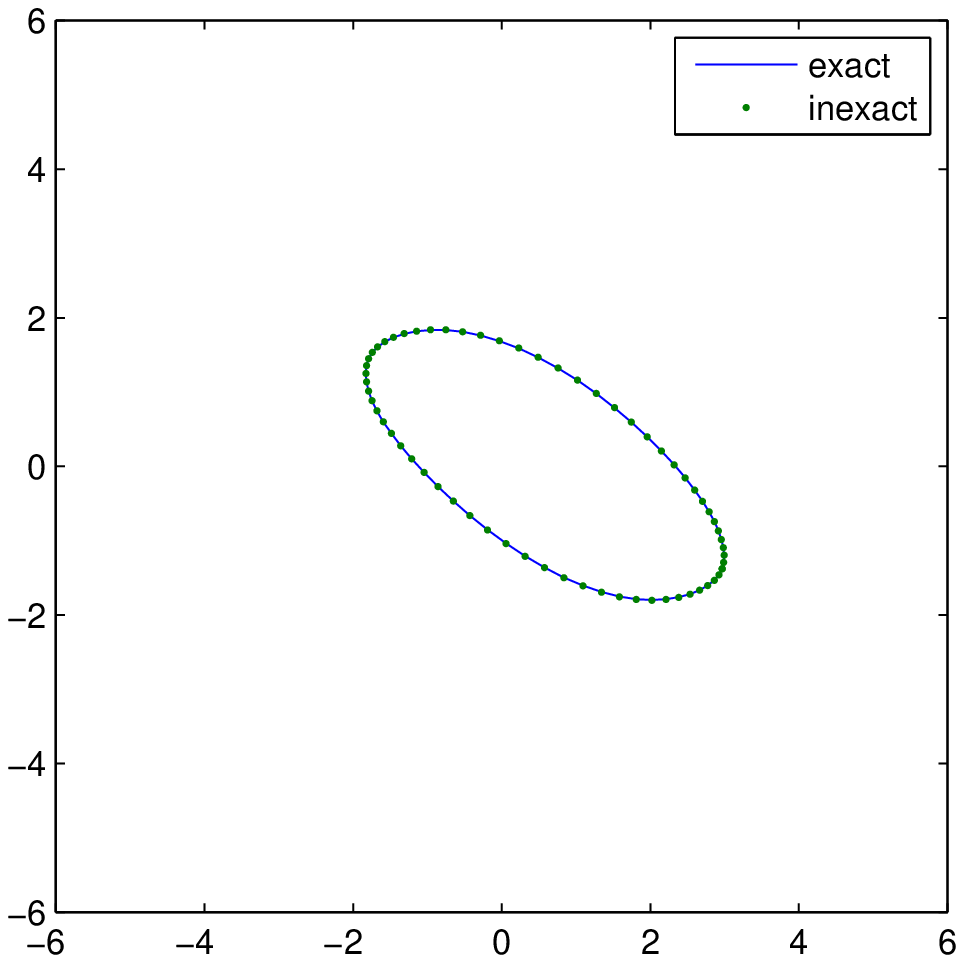}}
\subfigure[$t=0.8$]{\includegraphics[width = 2in]{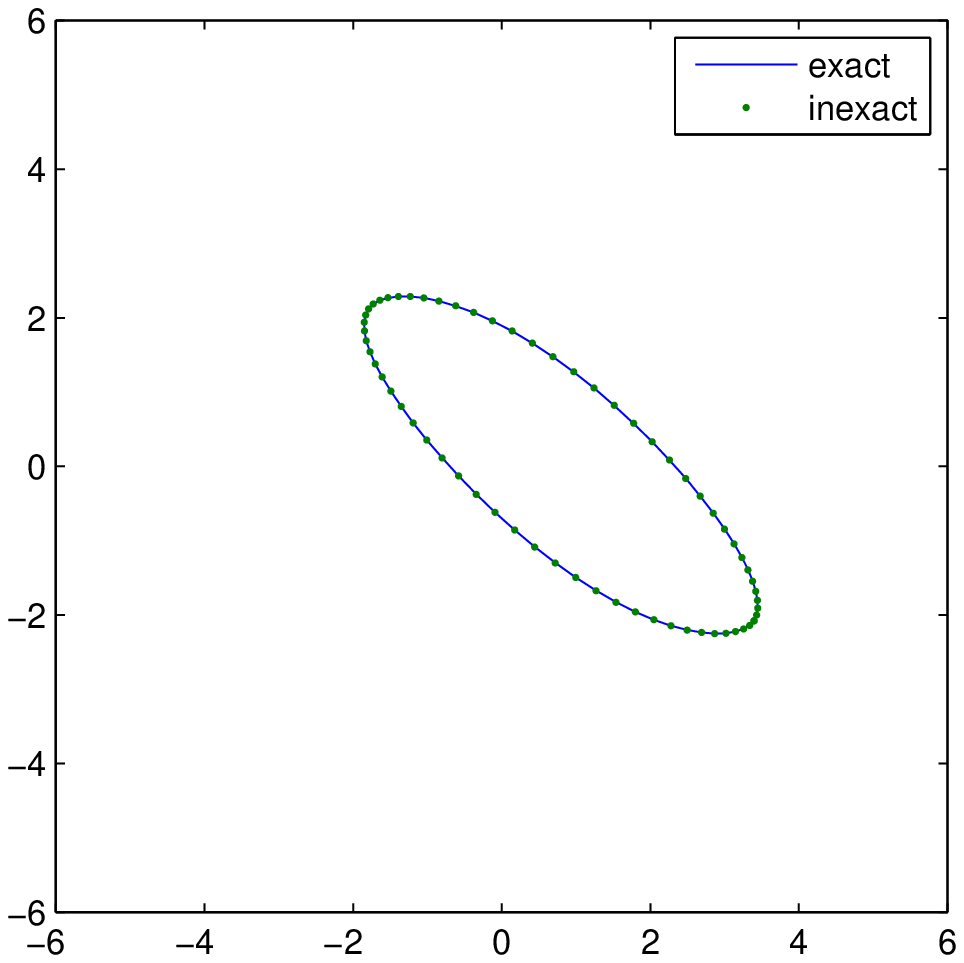}}
\subfigure[$t=1$]{\includegraphics[width = 2in]{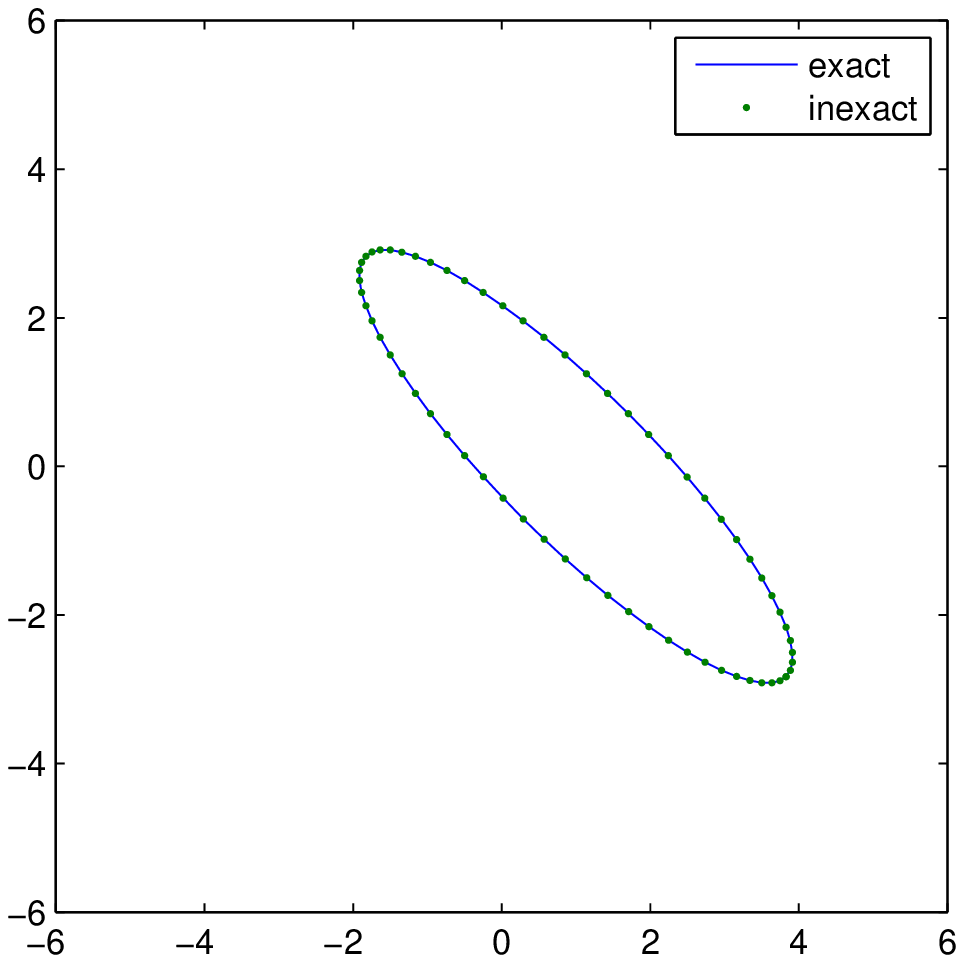}}
\caption{Comparison between exact and inexact matchings. $\sigma^2=0.1$ for the inexact matching.} 
\label{inexact1}
\end{figure}

\begin{figure}[h]
\subfigure[$t=0$]{\includegraphics[width = 2in]{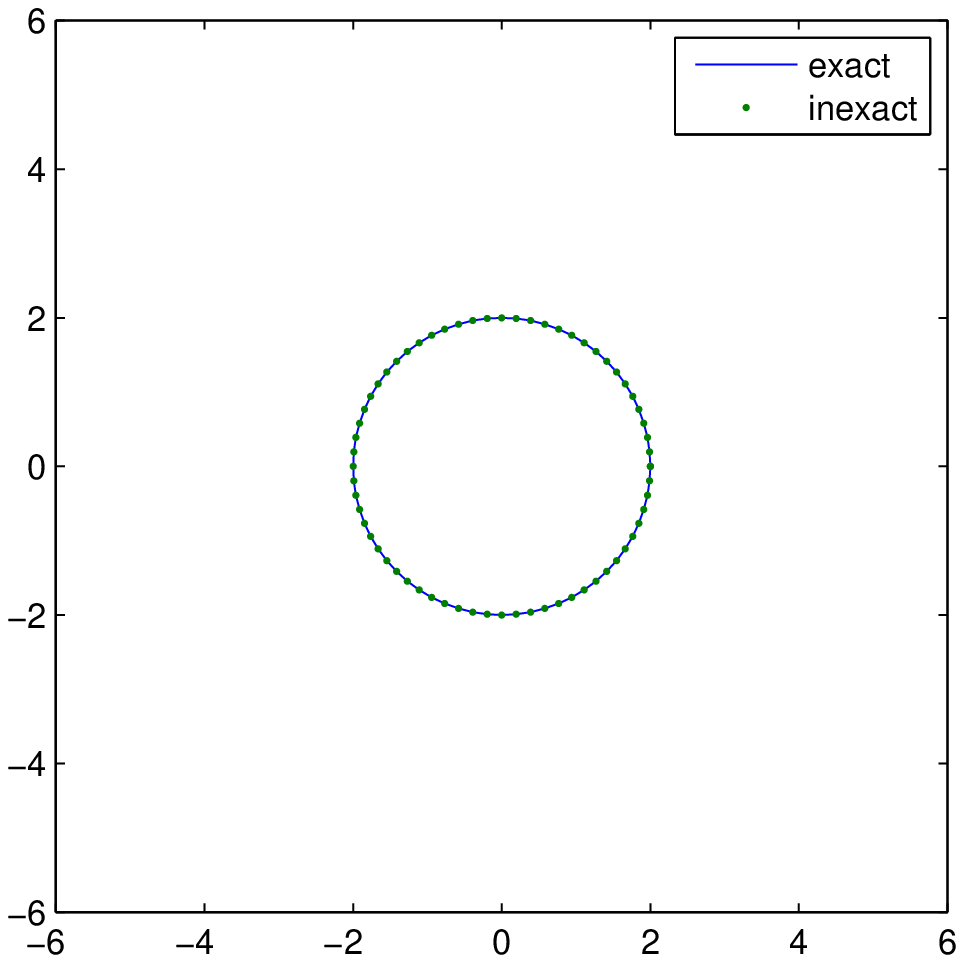}}  
\subfigure[$t=0.2$]{\includegraphics[width = 2in]{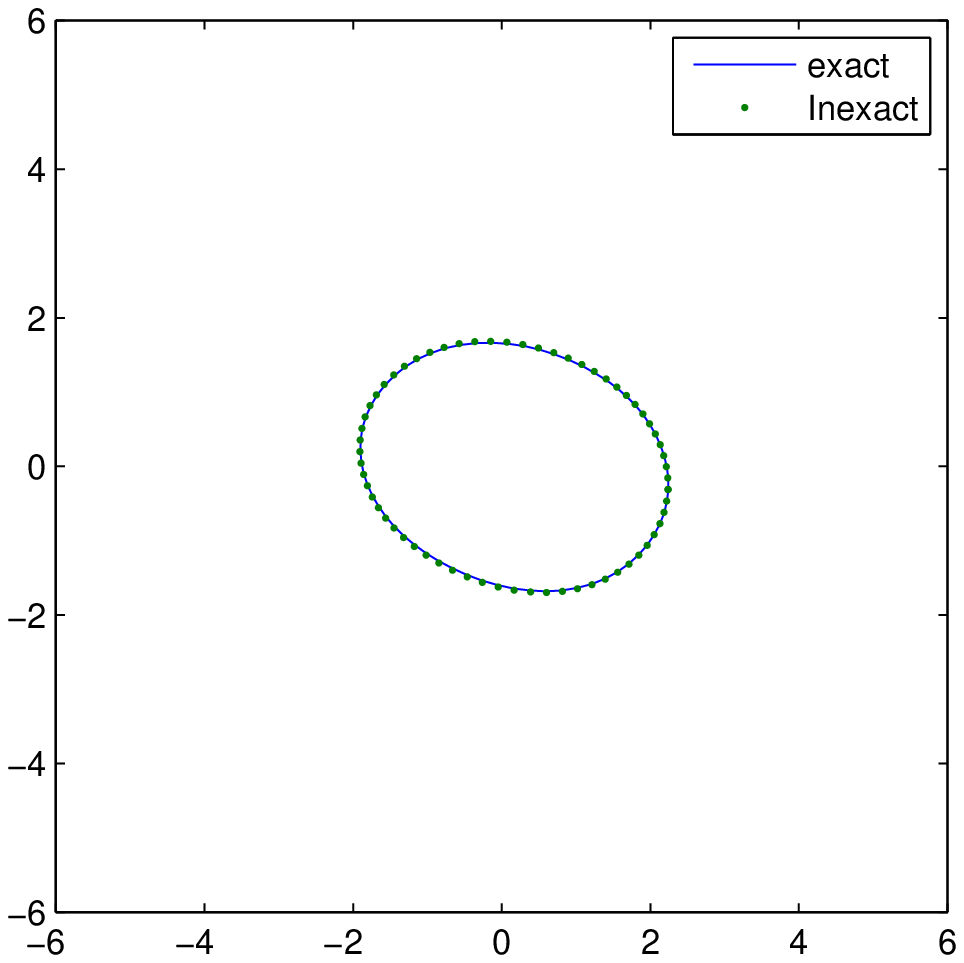}}
\subfigure[$t=0.4$]{\includegraphics[width = 2in]{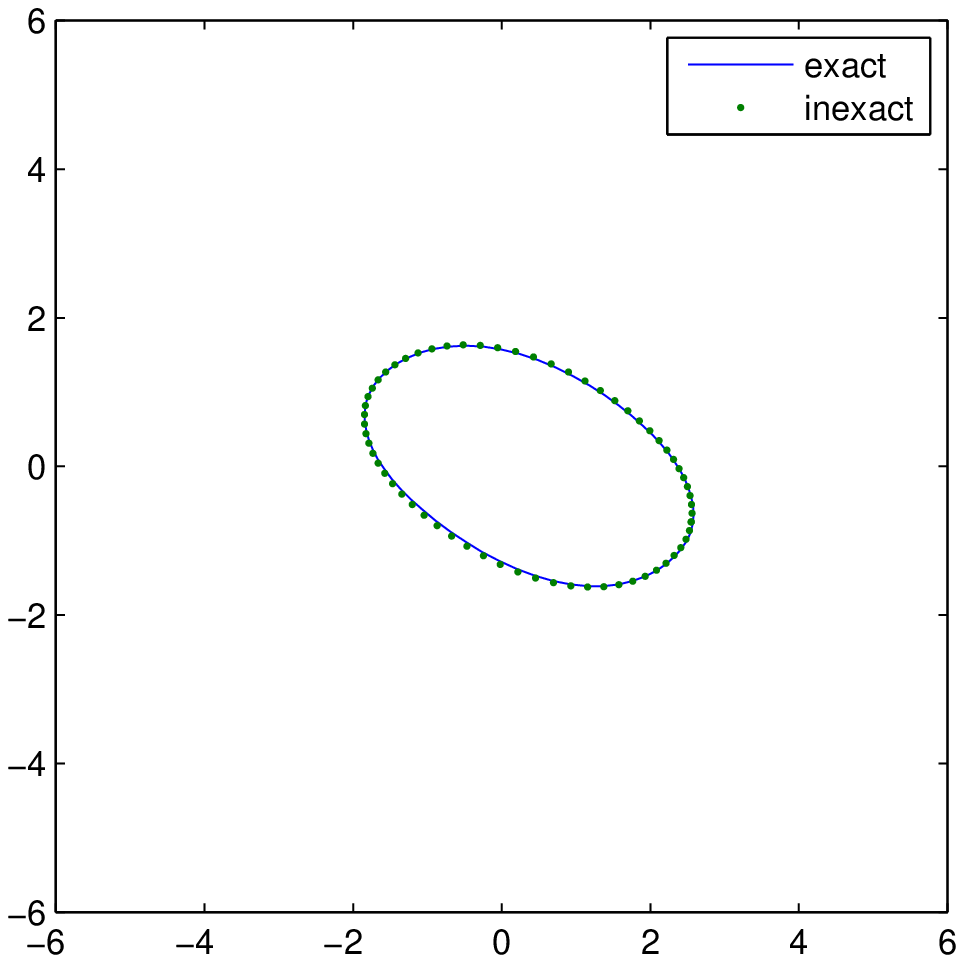}}\\  
\subfigure[$t=0.6$]{\includegraphics[width = 2in]{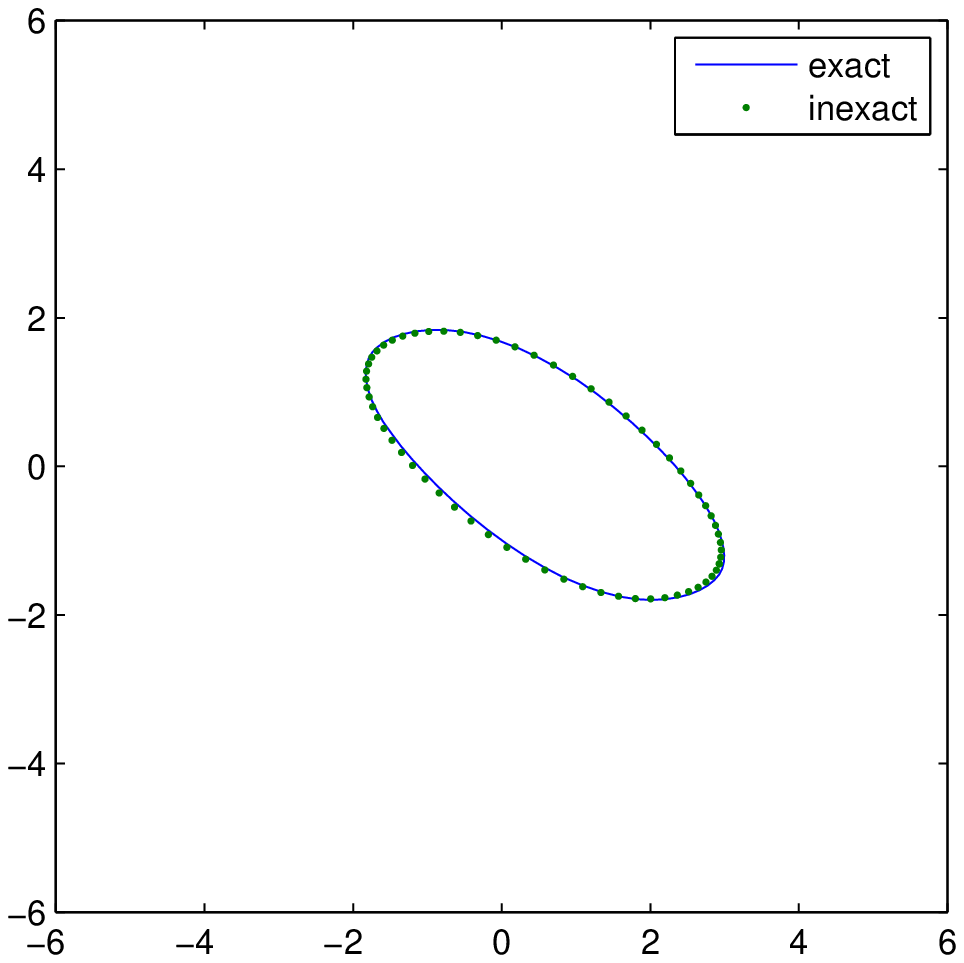}}
\subfigure[$t=0.8$]{\includegraphics[width = 2in]{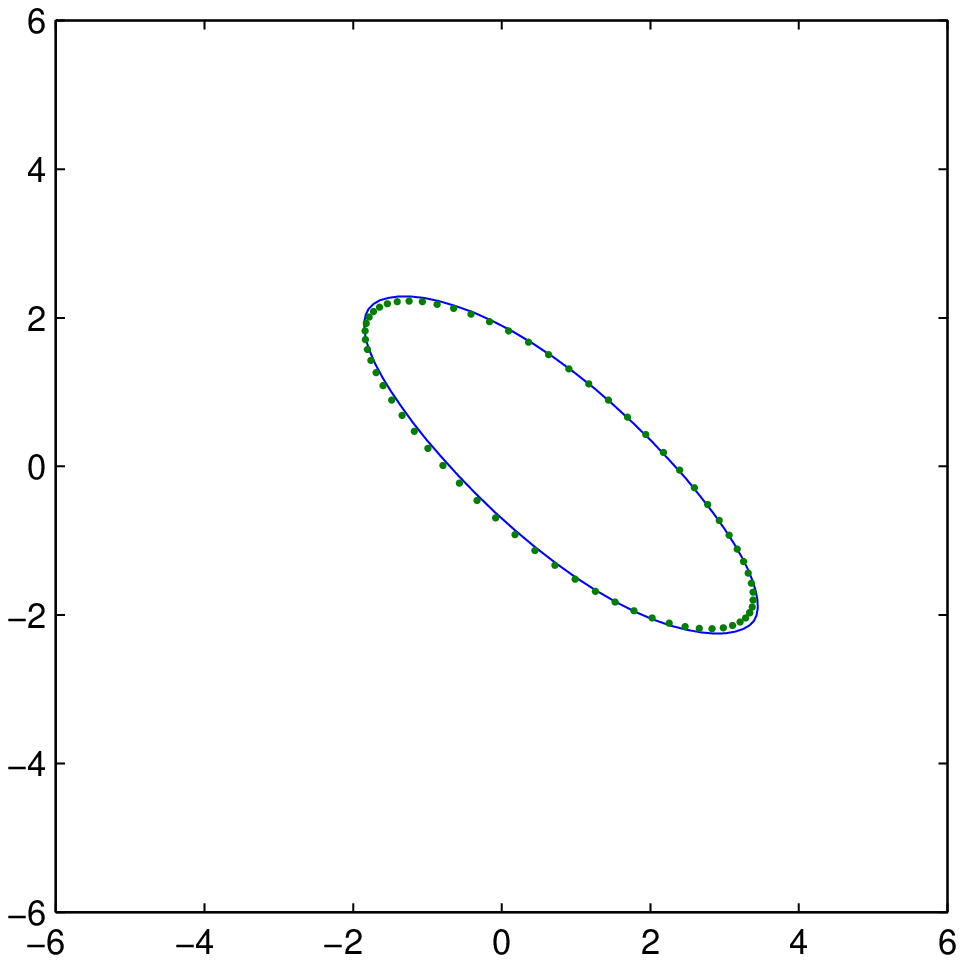}}
\subfigure[$t=1$]{\includegraphics[width = 2in]{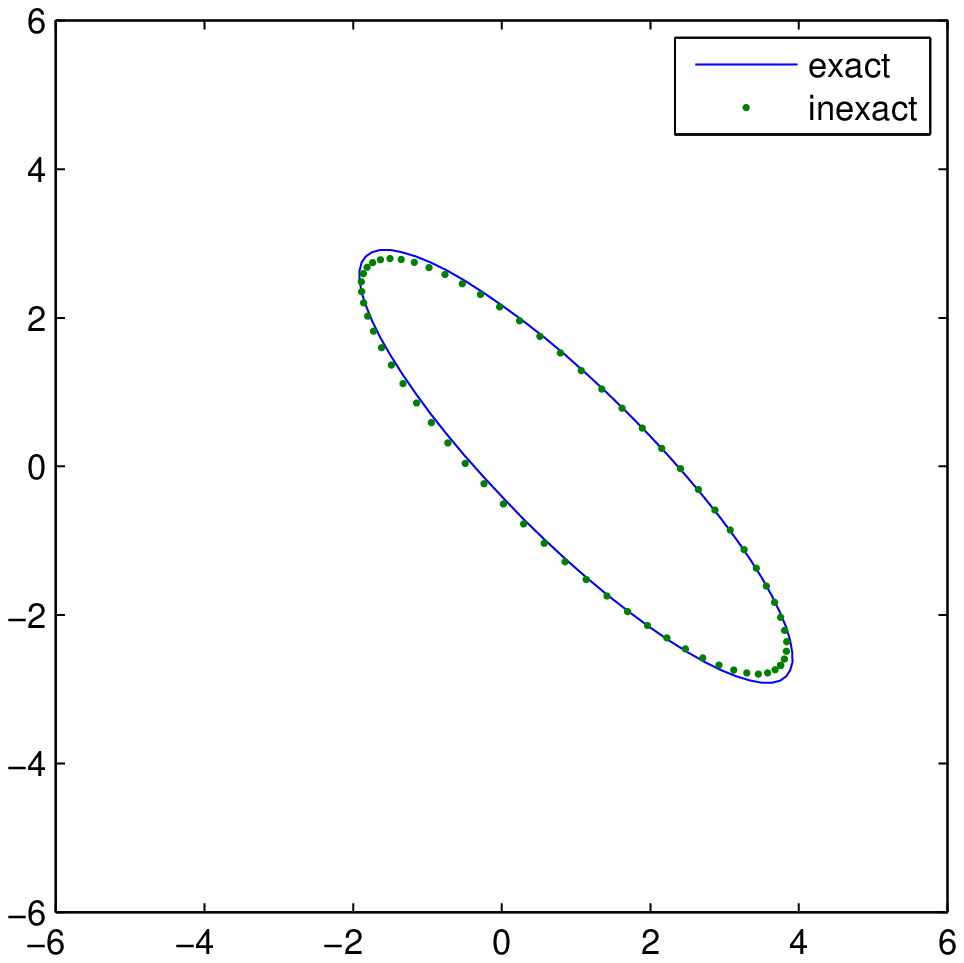}}
\caption{Comparison between exact and inexact matchings. $\sigma^2=0.5$ for the inexact matching.} 
\label{inexact2}
\end{figure}

We use $\sigma^2=0.1$ and $0.5$ for our experiments.  Table \ref{tab:ex-inex} shows the parameters used for the exact and the inexact matchings. The stopping criterion for the exact matching is $|I_1-I_1^{(k)}|_{\infty} < 10^{-3}$, while for the inexact matching is $|\bs{u}^{(k)}-\bs{u}^{((k-1))}|_{\infty} < 10^{-3}$. Figures \ref{inexact1} shows the comparison between the exact matching and the inexact matching with $\sigma^{2}=0.1$. Figure \ref{inexact2} is similar to Figure \ref{inexact1} but $\sigma^{2}=0.5$ for the inexact matching. We find that for $\sigma^2=0.1$, the inexact matching follows the exact matching closely, while for $\sigma^2=0.5$,  the two evolutions start showing discrepancy at $t=0.6$, and the difference is visible at $t=1$.  

\begin{table}[tbh]
\caption{Comparison between exact and inexact matchings}
\begin{tabular}{|l||c|c|c|}
\hline
\hline
& Exact  & Inexact ($\sigma^2=0.1$) & Inexact ($\sigma^2=0.5$)\\
\hline
Stopping criterion ($\epsilon$) &0.001 & 0.001 & 0.001\\ 
\hline
Hamiltonian $H$ &  46.5022   & 45.3927 & 41.3482\\
\hline
$|| I_1-I_1^{(k)}||_{\infty}$ & $9.5e^{-4}$ & 0.0042 & 0.0153\\
\hline
iteration number & 197 & 193 & 151\\
\hline
Search length $h$ & 0.3 & 0.4 & 0.1\\
\hline
\hline
\end{tabular}
\label{tab:ex-inex}
\end{table}

\section{Conclusion and Future Work}

We have presented a class of algorithms for template matching and its applications. The main difference between the present  algorithms and the traditional ones is that we treat the matching problem as an initial-value problem with unknown initial data. Other than solving the minimization problem for the Euler-Lagrangian equation, we iteratively solve the equivalent Euler-Poincar\'e equations with some initial guess for the initial data. The initial velocity that drives the diffeomorphism is obtained by an iterative process, manifesting a feedback-control loop. The advantages of the present algorithms include the use of the conical kernel in the particle system that limits the interaction between particles to accelerate the convergence and the availability of the implementation of fast-multipole method for solving the particle system. 

The extension of the current algorithms to solving three-dimensional template matching problems and the implementation of the fast-multiple methods for solving the particle system are currently under investigation. Other future work includes a different update scheme for algorithms under the stochastic setting and applying the algorithms for medical applications in two and three dimensions.

%

\section{Appendix I}\label{sec:app}
In this Appendix, we derive the expression of the optimal updating matrix $M$ described in Section \ref{sec:convergence}, but from the viewpoint of the estimation with least norm of covariance. The update of the velocity $\bs{u}$ is 
$${\bs u}_{i+1}={\bs u}_i+M\cdot(I_1-X_i)={\bs u}_i+M{\bs Q}_i.$$  Assuming that an unbiased guess for each step is given by:
$$E[{\bs Q}_i]=0, E({\bs u}_i)=\bs{u}.$$
The covariance matrix can be calculated as 
\begin{eqnarray}
Cov(\bs{u}_{i+1})&=&E[(\bs{u}_i-\bs{u}+M{\bs Q}_i)][(\bs{u}_i-\bs{u}+M{\bs Q}_i)^T]\\
            &=&Cov({\bs u}_i) + M\cdot Cov({\bs Q}_i)\cdot M^T +M\cdot Cov({\bs Q}_i,\bs{ u}_i) + Cov({\bs u}_i,{\bs Q}_i)\cdot M^T.
\end{eqnarray}
Then, we can find $M$ by calculating  
\begin{eqnarray}
&\frac{\partial\,\text{tr}(Cov({\bs u}_{i+1}))}{\partial M}=0;\\
&\Rightarrow 2[M\cdot Cov({\bs Q}_i) +Cov({\bs Q}_i,{\bs u}_i)]=0;\\
&\Rightarrow M=-Cov({\bs Q}_i,\bs{u}_i)\cdot Cov({\bs Q}_i)^{-1} \label{eq:M}.
\end{eqnarray}
${\bs Q}_i$ can be calculated as $${\bs Q}_i={\phi}(\bs{u})-{\phi}(\bs{u}_i)=-\frac{\partial {\phi}}{\partial {\bs u}}\cdot(\bs{u}_i-\bs{u})+O(|\bs{u}_i-\bs{u}|^2)\cdot 1_d.$$ If the second-order term is negligible compared to the first order term. We can then calculate
\begin{equation}\label{eq:cov1}
Cov({\bs Q}_i)=E({\bs Q_i}{\bs Q_i}^T)=\frac{\partial {\phi}}{\partial {\bs u}}\cdot Cov(\bs{u}_i)\cdot (\frac{\partial {\phi}}{\partial {\bs u}})^T,
\end{equation}
and
\begin{equation}\label{eq:cov2}
Cov (\bs {Q}_i, \bs{u}_i) = E[(\bs{u}-\bs{u}_i)({\bs Q}_i)^T] = -Cov(\bs{u}_i)(\frac{\partial \phi}{\partial \bs{u}})^T.
\end{equation}
Plugging equations (\ref{eq:cov1}) and (\ref{eq:cov2}) into equation (\ref{eq:M}) yields the first order estimation of the optimal $M$
$$M= (\frac{\partial {\phi}}{\partial\bs{u}})^{-1},$$
which is consistent with the local convergent analysis in Section \ref{sec:convergence}.

\end{document}